\documentclass[11pt,a4paper]{article}
\usepackage{amsmath, amssymb}
\usepackage{latexsym}
\usepackage{graphicx, epsfig}
\usepackage{color}

\usepackage{caption}
\usepackage{subcaption}

\newcommand{\SD}{\mathcal{S}_{\partial D}}
\newcommand{\SO}{\mathcal{S}_{\partial \Omega}}

\newcommand{\bef}{\begin{figure}}
	\newcommand{\enf}{\end{figure}}

\numberwithin{equation}{section}

\newtheorem{thm}{Theorem}[section]
\newtheorem{cor}[thm]{Corollary}

\def\nm{\noalign{\medskip}}

\newcommand{\qed}{\hfill \ensuremath{\square}}
\newcommand{\ds}{\displaystyle}
\newcommand{\pf}{\noindent {\sl Proof}. \ }
\newcommand{\p}{\partial}

\newcommand{\eqnref}[1]{(\ref {#1})}

\newcommand{\Cbb}{\mathbb{C}}

\newcommand{\Rbb}{\mathbb{R}}

\newcommand{\la}{\langle}
\newcommand{\ra}{\rangle}

\newcommand{\Kcal}{\mathcal{K}}

\newcommand{\Scal}{\mathcal{S}}

\newcommand{\Ga}{\alpha}
\newcommand{\Gb}{\beta}

\newcommand{\Ge}{\epsilon}
\newcommand{\Gve}{\varepsilon}

\newcommand{\Gvf}{\varphi}
\newcommand{\Gg}{\gamma}

\newcommand{\Gl}{\lambda}

\newcommand{\Gm}{\mu}
\newcommand{\Gv}{\nu}
\newcommand{\Gp}{\pi}
\newcommand{\Gt}{\theta}

\newcommand{\Gr}{\rho}

\newcommand{\Gs}{\sigma}

\newcommand{\Go}{\omega}

\newcommand{\Gz}{\zeta}
\newcommand{\GD}{\Delta}

\newcommand{\GG}{\Gamma}

\newcommand{\GO}{\Omega}

\newcommand{\beq}{\begin{equation}}
\newcommand{\eeq}{\end{equation}}
\newcommand{\ol}{\overline}

\begin{document}
\title{Neutral inclusions, weakly neutral inclusions, and an over-determined problem for confocal ellipsoids\thanks{\footnotesize This work was supported by NRF grants No. 2017R1A4A1014735 and 2019R1A2B5B01069967, JSPS KAKENHI grant No. 18H01126, NSF of China grant No. 11901523, and a grant from Central South University.}}

\author{Yong-Gwan Ji\thanks{Department of Mathematics and Institute of Applied Mathematics, Inha University, Incheon
22212, S. Korea. HK is also affiliated with Department of Mathematics and Statistics, Central South University, Changsha, Hunan, P.R. China. Permanent address of XL: College of Science, Zhejiang University of Technology, Hangzhou, 310023, P.R. China  (22151063@inha.edu, hbkang@inha.ac.kr, xiaofeili@zjut.edu.cn).} \and Hyeonbae Kang\footnotemark[2] \and Xiaofei Li\footnotemark[2] \and Shigeru Sakaguchi\thanks{Research Center for Pure and Applied Mathematics, Graduate School of Information Sciences, Tohoku University, Sendai, 980-8579, Japan (sigersak@tohoku.ac.jp).}}

\date{\today}
\maketitle

\begin{abstract}
An inclusion is said to be neutral to uniform fields if upon insertion into a homogenous medium with a uniform field it does not perturb the uniform field at all. It is said to be weakly neutral if it perturbs the uniform field mildly. Such inclusions are of interest in relation to invisibility cloaking and effective medium theory. There have been some attempts lately to construct or to show existence of such inclusions in the form of core-shell structure or a single inclusion with the imperfect bonding parameter attached to its boundary. The purpose of this paper is to review recent progress in such attempts. We also discuss about the over-determined problem for confocal ellipsoids which is closely related with the neutral inclusion, and its equivalent formulation in terms of Newtonian potentials. The main body of this paper consists of reviews on known results, but some new results are also included.
\end{abstract}



\section{Introduction}

This is a survey on recent progress in study on existence and construction of neutral and weakly neutral inclusions, and a related over-determined problem for confocal ellipsoids. The main body of the paper consists of reviews on known results with brief but coherent explanations. However, we include some new results as well.

To explain the problems related to the neutral and weakly neutral inclusion, let us consider the following conductivity problem:
$$
\mbox{(CP)} \
\begin{cases}
\nabla \cdot \Gs \nabla u = 0 \quad \mbox{in } \Rbb^d,\\
u(x)-a\cdot x= O(|x|^{-d+1}) \quad\mbox{as } |x| \to \infty,
\end{cases}
$$
for $d=2$ or $3$, where $a$ is a constant vector so that $-a = - \nabla (a\cdot x)$ is the background uniform field  and $\Gs$ is a piecewise constant function representing the conductivity distribution.

\begin{figure}[h!]
\begin{center}
\epsfig{figure=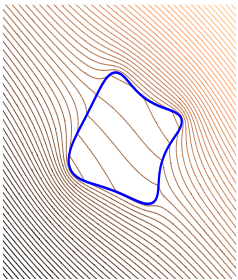,height=4.5cm,width=7cm}
\end{center}
\caption{Perturbation of the uniform fields; inside and outside the inclusion (with the boundary in blue)}\label{fig:perturbation}
\end{figure}

We first consider the problem (CP) when the conductivity distribution $\Gs$ is given by
\beq\label{Gssimply}
\Gs = \Gs_c \chi(D) + \Gs_m\chi(\Rbb^d \setminus D),
\eeq
where $D$ is a simply connected bounded domain in $\Rbb^d$ whose boundary $\p D$ is Lipchitz continuous. Here $\chi(D)$ denotes the characteristic function of $D$ ($\chi(\Rbb^d\setminus D)$ likewise), and $\Gs_c$ and $\Gs_m$ are constants representing conductivities of $D$ (the core) and $\Rbb^d\setminus D$ (the matrix), respectively.
In absence of the inclusion $D$, the solution to (CP) is nothing but $a\cdot x$. Thus, if we denote  by $u$ the solution to (CP) in presence of the inclusion, $u(x)- a\cdot x$ represents the perturbation occurred by insertion of the inclusion $D$ into the homogeneous medium with the uniform field $-a$. As we see from Fig. \ref{fig:perturbation}, the uniform field is perturbed outside (and inside) the inclusion.

It is known (see, e.g., \cite{AK07}) that the leading order term of the perturbation outside the inclusion can be expressed in terms of the dipolar expansion. In fact, we have the following expansion at infinity:
\beq\label{dipole}
u(x) - a \cdot x = \frac{1}{\Go_d} \frac{\la a , M x \ra}{|x|^d} +
 O(|x|^{-d}) \quad\mbox{as } |x| \to \infty,
\eeq
where $\Go_d$ is the surface area of the unit sphere in $\Rbb^d$ and $M=(M_{ij})$ is the $d \times d$ matrix determined by the domain $D$ and the conductivity contrast $\Gs_c/\Gs_m$. The matrix $M$ is called the polarization (or polarizability) tensor (PT in abbreviation, afterwards) associated with $D$. The PT is a signature of the existence of the inclusion $D$ and has been effectively used to detect some properties of the inclusion $D$, for which we refer to \cite{AK07}. It also plays an important role in the theory of composites and effective medium, for which we refer to \cite{AK07, mbook}.

If $D$ is a simply connected domain (or a union of simply connected domains), then $M$ is positive-definite if $\Gs_c - \Gs_m >0$ and negative-definite if $\Gs_c - \Gs_m <0$.  In fact, optimal bounds for PT, called the Hashin-Shtrikman bounds, are known, which will be explained in section \ref{sec:PT}.
Therefore, if $D$ is simply connected, then there is $\hat{x}=x/|x|$ such that $\la a , M \hat x \ra \neq 0$ and for such $x$ the following holds
\beq\label{dipole2}
|u(x) - a \cdot x| \ge C |x|^{-d+1} \quad\mbox{as } |x| \to \infty
\eeq
for some $C>0$. The dipolar expansion \eqnref{dipole} shows that in general the solution $u$ to (CP) admits the following:
\beq\label{general}
u(x)-a\cdot x=O(|x|^{-d+1}) \quad\mbox{as } |x|\to\infty.
\eeq
Furthermore, \eqnref{dipole2} shows that the decay rate $O(|x|^{-d+1})$ cannot be replaced by a faster rate, say $O(|x|^{-d})$.

\begin{figure}[t!]
\begin{center}
\epsfig{figure=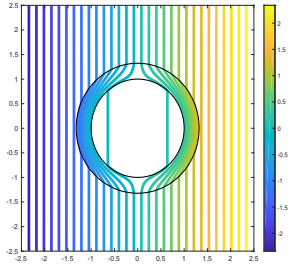, height=4cm, width=5cm}
\end{center}
\caption{Neutral inclusion: The uniform field is not perturbed}\label{fig:Hashine}
\end{figure}

However, if the inclusion is of a core-shell structure, then the situation can be quite different. Let $D$ be a bounded domain and $\GO$ be a bounded domain containing $\ol{D}$ so that $(D, \GO)$ becomes a coated structure or a core-shell structure. Suppose that the conductivity distribution is given by
\beq\label{Gsdouble}
\Gs = \Gs_c \chi(D) + \Gs_s \chi(\GO \setminus D) + \Gs_m\chi(\Rbb^d \setminus \GO).
\eeq
In particular, if $(D, \GO)$ is a pair of concentric disks or balls, and if the conductivities $\Gs_c$, $\Gs_s$ and $\Gs_m$ are scalars and satisfy
\beq\label{effective}
(d-1 + \Gs_c/\Gs_s)(\Gs_m/\Gs_s - 1) + f (1 - \Gs_c/\Gs_s)(\Gs_m/\Gs_s + d-1) =0
\eeq
for $d=2$ or $3$, where $f=|D|/|\GO|$ (the volume fraction), then the solution to (CP) satisfies
\beq\label{neutral}
u(x)-a\cdot x \equiv 0 \quad\mbox{in } \Rbb^d \setminus \GO,
\eeq
namely, the uniform field is not perturbed at all (see Fig. \ref{fig:Hashine}).
In fact, with the conductivity given by \eqnref{Gsdouble}, the solution $u$ to (CP) is harmonic in $\Rbb^d \setminus (\p D \cup \p \GO)$, and along the interfaces $\p D$ and $\p \GO$ it satisfies the transmission conditions: continuity of the potential and continuity of the flux, namely,
\beq\label{interface}
\Gs_c \p_\nu u|_- = \Gs_s \p_\nu u|_+ \mbox{ on } \p D, \quad \Gs_s \p_\nu u|_- = \Gs_m \p_\nu u|_+ \mbox{ on } \p \GO,
\eeq
where the subscripts $+$ and $-$ indicate the limits from outside and inside $D$ (or $\GO$), respectively. If $D$ and $\GO$ are concentric disks (or balls), one can use spherical harmonics to find the solution explicitly to satisfy these interface conditions, and show that \eqnref{effective} implies \eqnref{neutral}.

This easy-to-prove fact was first discovered by Hashin \cite{H1}, and significance of the discovery lies in its implications.
Since insertion of inclusions does not perturb the outside uniform field, the effective conductivity of the assemblage filled with such inclusions of many different scales is the same as $\Gs_m$ (the conductivity of the matrix) satisfying \eqnref{effective}. It is also proved that such an effective conductivity is one of the Hashin-Shtrikman bounds on the effective conductivity of arbitrary two-phase composites \cite{H1, HS} (see also \cite{mbook}).

The inclusion $(D, \GO)$ of core-shell structure (or any other structure), which does not perturb the uniform field $-a$ upon its insertion, that is, satisfying \eqnref{neutral}, is said to be neutral to the field $-a$. If the inclusion is neutral to all uniform fields, it is said to be neutral to multiple uniform fields. The concentric disks (or balls) satisfying \eqnref{effective} is neutral to multiple uniform fields. If $\Gs_m$ is anisotropic ($\Gs_c$ and $\Gs_s$ are scalars), then confocal ellipsoids whose common foci are determined by $\Gs_m$ can be neutral to multiple fields. We include a proof of this fact in section \ref{sec:neutral}.

Now the question is if there are coated inclusions other than concentric disks or balls neutral to multiple uniform fields (or confocal ellipses or ellipsoids when $\Gs_m$ is anisotropic), more generally, if we can coat a given domain $D$ of general shape by a shell so that the resulting inclusion $(D, \GO)$ is neutral to multiple uniform fields. The answer is proven to be no in two dimensions. In fact, it has been proved that if $(D, \GO)$ is neutral to multiple uniform fields, then $\GO$ and $D$ are concentric disks if $\Gs_m$ is isotropic, and
confocal ellipses if $\Gs_m$ is anisotropic (and the foci are determined by $\Gs_m$). This was proved by Milton-Serkov \cite{MS} when $\Gs_c=0$ or $\infty$, and by Kang-Lee \cite{KL14} when $\Gs_c$ is finite. Since these two dimensional results are proved using either the Riemann mapping or existence of harmonic conjugates, the methods of proofs cannot be extended to three dimensions. It is worth mentioning that there are many different shapes of coated inclusions neutral to a single uniform field as shown in two dimensions in \cite{JM, MS}. In three dimensions, it is proved in \cite{KLS} that the coated inclusion $(D, \GO)$ being neutral to multiple fields is equivalent to existence of a solution to a certain over-determined problem defined on $\GO \setminus \ol{D}$. It is then proved as a consequence that if $\Gs_m$ is isotropic, then the only inclusions of core-shell structure is a pair of concentric balls. Extension of this result to the anisotropic case has not been proved and is open. We will review on neutral inclusions and the related over-determined problem in section \ref{sec:neutral}. We also include in the same section a proof of their equivalence to a certain formulation in terms of Newtonian potentials.

Other than applications to the theory of composite as explained earlier, there is another interest in neutral inclusions in relation to invisibility cloaking. The neutrality condition \eqnref{neutral} means that the uniform field is unperturbed at all outside the inclusion, namely, there is no difference of the field with or without the inclusion. It means that the inclusion is invisible from the probe by uniform fields. This was also observed in \cite{kerker}.  Recently, the idea of neutrally coated inclusions has been extended to construct multi-coated circular structures which are neutral not only to uniform fields but also to fields of higher order up to $N$ for a given integer $N$ \cite{AKLL1}. It was proved there that the multi-coated structure combined with a transformation optics dramatically enhances the near cloaking of \cite{kohn1}.

Since there is no coated inclusion neutral to multiple fields (invisible by uniform fields), we may ask if there are inclusions which are vaguely visible by uniform fields. They are weakly neutral inclusions.
In general, the solution to (CP) satisfies the decay condition $u(x)-a\cdot x=O(|x|^{1-d})$ at $\infty$ and if the inclusion is neutral, then $u(x)-a\cdot x \equiv 0$ outside the inclusion. This property means that the field outside the inclusion is not perturbed even though the inclusion is inserted. The weakly neutral inclusions perturb the fields mildly:
\beq\label{wneutral}
u(x)-a\cdot x=O(|x|^{-d}) \quad\mbox{as } |x|\to\infty.
\eeq
If \eqnref{wneutral} holds for all constant vectors $a$, then the inclusion is said to be weakly neutral to multiple uniform fields.
According to \eqnref{dipole}, in order for \eqnref{wneutral} to hold for all $a$, the corresponding PT must vanish. Thus weakly neutral inclusions are PT cancelling ones. We now formulate the weakly neutral inclusion problem:

\begin{itemize}
\item[]
{\bf Weakly Neutral Inclusion Problem}. {\sl Given a domain $D$ can we find a domain $\GO$ containing $\ol{D}$ so that the resulting inclusion of core-shell structure becomes weakly neutral to multiple uniform fields, or equivalently, its PT vanishes. }
\end{itemize}

In section \ref{sec:weak}, we present two classes of domains which admit coatings so that the resulting inclusions of core-shell structure become weakly neutral to multiple uniform fields. One class is the collection of domains $D$ such that the coefficients $b_D$ vanish. Here $b_D$ is the leading coefficient of the conformal mapping from the exterior of the unit disk onto the exterior of $D$ (see \eqnref{f}). For such domains we construct the coating explicitly. This is a new result. The other class is that of small perturbations of disks, for which it is proved in \cite{KLS2D} that there are coatings so that the resulting inclusions become weakly neutral to multiple uniform fields.

There is yet another way, other than coating, to achieve weak neutrality. It is by introducing imperfect bonding parameter on the boundary of the given domain. We review the result of \cite{KL19} on this in section \ref{sec:imperfect}.

This paper is organized in the following way. Section \ref{sec:PT} is to review general properties of the PT, including the Hashin-Shtrikman bounds. In section \ref{sec:neutral}, we discuss about neutral inclusions and related over-determined problem for confocal ellipsoids and an equivalent formulation in terms of the Newtonian potential. We also include a discussion on Neumann ovaloids. Section \ref{sec:weak} is to discuss on the weakly neutral inclusion problem. Section \ref{sec:imperfect} is for discussion on the construction of weakly neutral inclusion by imperfect bonding parameters.

\section{Layer potentials and polarization tensors}\label{sec:PT}

In this section we represent the PT appearing in the dipolar expansion \eqnref{dipole} in terms of layer potentials and recall the optimal Hashin-Shtrikman (HS) bounds on traces of the PT and its inverse.

\subsection{Layer potentials}\label{subsec:layer}

Let $\GG(x)$ be the fundamental solution to the Laplacian, that is, $\GG(x)= 1/(2\pi) \log |x|$ in two dimensions, and $\GG(x)=- (4\pi |x|)^{-1}$ in three dimensions. Let $D$ be a bounded simply connected domain with a Lipschitz continuous boundary. The single layer potential of a function $\Gvf\in H^{-1/2}(\p D)$ (the $L^2$ Sobolev space of order $-1/2$ on $\p D$) is defined by
\beq
\Scal_{\p D} [\Gvf](x) := \int_{\p D} \GG(x-y) \Gvf(y)\, dS(y), \quad x\in \Rbb^d,
\eeq
where $dS$ is the line or surface element on $\p D$. Let $\p_{\nu}$ denote the outward normal derivative on $\p D$. It is well known (see, for example, \cite{AK07}) that the following jump relation holds:
\beq\label{jump}
\p_\nu \Scal_{\p D}[\Gvf](x)\big|_{\pm} = \left( \pm \ds\frac{1}{2}I+\Kcal_{\p D}^{*} \right)[\Gvf](x), \quad\mbox{a.e. } x\in\p D,
\eeq
where $I$ is the identity operator and $\Kcal_{\p D}^{*}$ is the operator defined by
\begin{equation*}
\Kcal_{\p D}^{*}[\Gvf](x)= \frac{1}{\Go_d} \int_{\p D} \frac{\la x-y,\nu(x)\ra}{|x-y|^d} \Gvf(y)\, dS(y).
\end{equation*}
Here, $\la \, , \, \ra$ the scalar product in $\Rbb^d$. The boundary integral operator $\Kcal_{\p D}^{*}$ is called the Neumann-Poincar\'e (NP) operator.

\subsection{Polarization tensors}\label{subsec:PT}

Let $u_l$, $1 \le l \le d$, be the solution to (CP) when $a \cdot x= x_l$ and the conductivity distribution $\Gs$ is given by \eqnref{Gssimply}. Then it is known (see, e.g., \cite{AK07}) that $u_l$ can be represented as
\beq\label{repsimple}
u_l(x)=x_l+\Scal_{\p D} [\Gvf^{(l)}](x),\quad x\in\Rbb^d,
\eeq
where $\Gvf^{(l)}$ is the unique solution in $H^{-1/2}_0(\p D)$ ($H^{-1/2}(\p D)$ functions with the mean zero) to the integral equation
\beq
\left( \frac{\Gs_c+\Gs_m}{2(\Gs_c-\Gs_m)}I - \Kcal_{\p D}^* \right) [\Gvf^{(l)}] = \nu_l,
\eeq
where $\nu_l$ is the $l$-th component of the outward unit normal vector field $\nu$ on $\p D$. By expanding out the term $\Scal_{\p D} [\Gvf^{(l)}](x)$ in \eqnref{repsimple} as $|x| \to \infty$, we see that the PT $M=M(D)=(m_{ll'})_{l,l'=1}^d$ in this case is given by
\beq\label{PTsimple}
m_{ll'}=\int_{\p D} x_{l'}\Gvf^{(l)} \, dS.
\eeq

If the conductivity distribution $\Gs$ is given by \eqnref{Gsdouble}, the solution $u_l$ can be represented as
$$
u_l(x)=x_l+\SD[\Gvf_1^{(l)}](x)+\SO[\Gvf_2^{(l)}](x),\quad x\in\Rbb^d,
$$
where $(\Gvf_1^{(l)},\Gvf_2^{(l)})\in H^{-1/2}_0(\p D)\times  H^{-1/2}_0(\p \GO)$ is the unique solution to the system of integral equations
\begin{equation}\label{pp}
\begin{bmatrix}
-\Gl I+\Kcal_{\p D}^* & \p_\nu\Scal_{\p\GO} \\
\p_\nu\Scal_{\p D} & -\mu I+\Kcal_{\p\GO}^*
\end{bmatrix}
\begin{bmatrix}
\Gvf_1^{(l)}\\
\Gvf_2^{(l)}
\end{bmatrix}
=-\begin{bmatrix}
\nu_l^{\p D}\\
\nu_l^{\p \GO}
\end{bmatrix}.
\end{equation}
Here we denote the unit outward normal vector fields on $\p D$ and $\p \GO$ by $\nu^{\p D}$ and $\nu^{\p\GO}$, respectively. The numbers $\Gl$ and $\Gm$ are given by
\beq\label{lambdamu}
\Gl=\frac{\Gs_c+\Gs_s}{2(\Gs_c-\Gs_s)}\quad\mbox{and}\quad \mu=\frac{\Gs_s+\Gs_m}{2(\Gs_s-\Gs_m)}.
\eeq
For unique solvability of the integral equation we refer to \cite{KLS2D}. In this case, the PT $M=M(D,\GO)=(m_{ll'})_{l,l'=1}^d$ of the core-shell structure $(D,\GO)$ is given by
\beq\label{PT0}
m_{ll'}=\int_{\p D}x_{l'}\Gvf_1^{(l)} \, dS+\int_{\p \GO}x_{l'} \Gvf_2^{(l)} \, dS.
\eeq
It is known that $M$ is a symmetric matrix (see, e.g., \cite{AK07}).

\subsection{Hashin-Shtrikman bounds}\label{subsec:HS}

If the conductivity distribution is given by \eqnref{Gssimply}, then the following optimal bounds on traces of the PT $M$ and its inverse hold: with $k = \Gs_c/\Gs_m$,
\beq\label{HSupper}
 \mbox{Tr}(M) < |D| (k-1) (d-1 + \frac{1}{k}),
\eeq
and
\beq\label{HSlower}
 |D| \mbox{Tr}(M^{-1})\leq\frac{d-1+k}{k-1},
\eeq
where Tr stands for trace. These bounds are obtained by Lipton \cite{Lipton93} under the assumption of periodicity, and by Capdeboscq-Vogelius \cite{CV05} without assumption of periodicity, and called the Hashin-Shtrikman (HS) bounds after the names of scientists who found optimal bounds of effective properties of two phase composites, as described in the paragraph right after \eqnref{interface}.

\begin{figure}[t!]
\begin{center}
\epsfig{figure=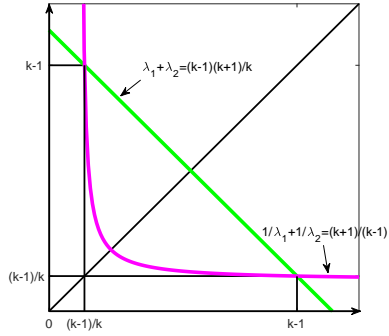, height=5cm}
\end{center}
\caption{Hashin-Shtrikman bounds for the PT}\label{HSbound}
\end{figure}

The first one is an upper bound (the green line in Fig. \ref{HSbound}) and the second one is a lower one (the pink curve in Fig. \ref{HSbound}).  The upper bound is never attained by a domain, while the lower bound is attained by ellipses and ellipsoids, and the converse is also true. In fact, it is proved in \cite{KM06, KM08} that the simply  connected domain whose PT satisfies the lower HS-bounds is an ellipse in two dimensions and an ellipsoid in three dimensions. This is an isoperimetric inequality for PT and a generalized version of the P\'olya-Szeg\"o conjecture \cite{PS51}. The original P\'olya-Szeg\"o conjecture asserts that the PT attains its minimal trace on and only on disks or balls. The constant trace lines of the PT are those parallel to the green line in Fig. \ref{HSbound}. Thus the minimal trace is attained at the point of tangency of the line parallel to the green line to the pink curve. The generalized version asserts that if eigenvalues of the PT lies on the pink curve, then the domain must be an ellipse or an ellipsoid. The original P\'olya-Szeg\"o conjecture is now proved as a simple consequence of the generalized version. See Theorem \ref{thm:ell} of this paper for more discussion on this. The bounds \eqnref{HSupper} and \eqnref{HSlower} are optimal in the sense that any matrix satisfying \eqnref{HSupper} and \eqnref{HSlower} is actually the PT associated with a domain (see \cite{ACKKL, CV03} for proofs).

\section{Neutral inclusions and an over-determined problem}\label{sec:neutral}

In this section, the conductivity distribution $\Gs$ is given by \eqnref{Gsdouble} with the inclusion $(D,\GO)$ of core-shell structure. We assume that the conductivity of the matrix, $\Gs_m$, is in general anisotropic, i.e., a symmetric matrix. We review the result saying that if $\Gs_m$ is isotropic, i.e., its eigenvalues are all the same, then the only inclusion of the core-shell structure neutral to multiple uniform fields is concentric balls. We also prove the equivalence of neutral inclusion problems with an over-determined problem for confocal ellipsoids, and an equivalent formulation of the problem using the Newtonian potentials. In relation to these problems, we include a subsection on quadrature domains and Neumann ovaloids.

\subsection{An over-determined problem for confocal ellipsoids}\label{subsec:confocal}

It was proved in \cite{KLS} that if $(D,\GO)$ is neutral to multiple uniform fields and $\Gs_c>\Gs_s$, then the following over-determined problem admits a solution:
$$
\mbox{(ODP)} \
\left\{
\begin{array}{ll}
\ds \GD w= 1 \quad &\mbox{in } \GO \setminus \overline{D}, \\
\nm
\ds \nabla w = 0 \quad &\mbox{on } \p \GO , \\
\nm
\ds \nabla w(x)= Ax + b \quad &\mbox{on } \p D,
\end{array}
\right.
$$
where $A$ is a symmetric matrix and $b$ is a constant vector, provided that $\p D$ is connected and $\Rbb^3 \setminus \overline{D}$ is simply connected. This problem is over-determined since $\nabla w$ is prescribed on $\p\GO$ and $\p D$. The matrix $A$ is determined by $\Gs_m$. If $\Gs_m$ is isotropic for example, so is $A$.

Let us briefly recall the proof. Suppose, after diagonalization, that
\beq\label{Gsmaniso}
\Gs_m = \mbox{diag} [\Gs_{m,1}, \Gs_{m,2}, \Gs_{m,3}].
\eeq
Let $e_j$, $j=1,2,3$, be the standard basis of $\Rbb^3$ and let $u_j$ be the solution to (CP) when $a=e_j$.
The inclusion $(D,\GO)$ being neutral to multiple uniform fields means that $u_j(x)-x_j=0$ in $\Rbb^3 \setminus \GO$ for $j=1,2,3$.
Let
\beq \label{beta}
\Gb_j := \frac{\Gs_{m,j}}{\Gs_s}-1,
\eeq
and
\beq\label{wj}
v =  (\Gb_1^{-1} u_1, \Gb_2^{-1} u_2, \Gb_3^{-1} u_3)^T.
\eeq
The crux of the proof in \cite{KLS} lies in proving that $v$ is linear inside $D$. In fact, it is proved that $v(x)=c_0x+b_0$ ($x \in D$) for some constant $c_0$ and vector $b_0$. It is here where the assumption $\Gs_c > \Gs_s$ is needed\footnote{We believe it is true without this assumption even though we do not know how to prove it.}. It is then shown that $\nabla v$ is symmetric, and hence, thanks to the assumption that $\p D$ is connected and $\Rbb^3 \setminus \overline{D}$ is simply connected, there is a function $\psi$ in $\overline{\GO}\setminus D$ such that $v=\nabla \psi$. Moreover, $\Delta \psi = \sum_{j=1}^3 \Gb_j^{-1} +1$ in $\GO \setminus \overline{D}$. Thus $w$, defined by
\beq\label{definition_of_w}
w(x):= \psi(x)  - \frac{1}{2} \sum_{j=1}^3 \Gb_j^{-1}x_j^2 ,
\eeq
is the solution to (ODP) with $A=c_0I - \mbox{diag}[\Gb_1^{-1}, \Gb_2^{-1}, \Gb_3^{-1}]$ ($I$ is the identity matrix). If $\Gs_m$ is isotropic, so is $A$ as mentioned before. The converse is also true, namely, if (ODP) admits a solution, then $(D,\GO)$ is neutral. For this, see Theorem \ref{thm:equiv} below.

\medskip

\noindent {\bf Remark} {\it
The assumption that $\p D$ is connected and $\mathbb R^3 \setminus \overline{D}$ is simply connected in {\rm \cite[Theorem 1.2]{KLS}} can be replaced with the weaker one that $\Omega \setminus \overline{D}$ is connected.
Indeed, instead of using Stokes' theorem, by combining the formula  {\rm \cite[(2.18)]{KLS}} with the fact that $v(x) = c_0x + b_0$ $(x \in D)$, we see that the function $\psi$ is explicitly given by
\beq\label{without Stokes}
\psi(x) =  c_0\left(1-\frac {\sigma_c}{\sigma_s}\right)\int_D\Gamma(x-y) dy + \frac 12 x\cdot Bx + \int_\Omega\Gamma(x-y) dy.
\eeq
Hence the function $w = w(x)$ is given by
\beq\label{explicit representation of w}
 w(x) = c_0\left(1-\frac {\sigma_c}{\sigma_s}\right)\int_D\Gamma(x-y) dy + \int_\Omega\Gamma(x-y) dy.
\eeq
For a domain $D$ in three dimensions and a domain $\Omega$ containing $\overline{D}$, the assumption that $\Omega \setminus \overline{D}$ is connected is really more general than that $\p D$ is connected and $\mathbb R^3 \setminus \overline{D}$ is simply connected.  In fact,  this general assumption allows us to choose the genus of a closed surface $\p D$ arbitrarily.  If the genus does not equal zero,  $\mathbb R^3 \setminus \overline{D}$ is not simply connected,  but $\Omega \setminus \overline{D}$ is connected.
 }

\medskip

Note that if $\GO$ and $D$ are concentric balls centered at the origin whose respective radii are $r_e$ and $r_i$, then the solution $w$ to (ODP) is given by

\beq
w(x) = \frac{r_e^3}{3|x|} + \frac{1}{6} |x|^2.
\eeq
In this case, $b=0$ and $A=\frac{1}{3}(-r_e^3/r_i^3+1)I$ which is isotropic. We emphasize that $w$ is radial in this case.

It is shown in \cite{KLS} that confocal ellipsoids admit a solution to (ODP). In fact,
if $\p D$ is an ellipsoid given by
\beq\label{ellD}
\frac{x^2_1}{c^2_1}+\frac{x^2_2}{c^2_2}+\frac{x^2_3}{c^2_3}=1,
\eeq
the confocal ellipsoidal coordinate $\Gr$ is given by
\beq\label{ellGO}
\frac{x^2_1}{c^2_1+\rho}+\frac{x^2_2}{c^2_2+\rho}+\frac{x^2_3}{c^2_3+\rho}=1,
\eeq
and the confocal ellipsoid $\p\GO$ is given by $\rho=\rho_0$ for some $\Gr_0>0$.
Let
\beq\label{gGr}
g(\rho) = (c_1^2+\rho)(c_2^2+\rho)(c_3^2+\rho),
\eeq
and
\beq\label{Gvfj}
\Gvf_j (\rho)= \int_\rho^\infty \frac{1}{(c_j^2+s)\sqrt{g(s)}} ds, \quad j=1,2,3.
\eeq
Then the function $w$, defined by
\beq
w(x) = \frac{1}{2} \int_\rho^\infty \frac{1}{\sqrt{g(s)}} ds - \frac{1}{2} \sum_{j=1}^3 \Gvf_j(\rho) x_j^2 + \frac{1}{2} \sum_{j=1}^3 \Gvf_j(\rho_0) x_j^2,
\eeq
is a solution of (ODP) with
$$
A= \mbox{diag} [ \Gvf_1(\rho_0) - \Gvf_1(0), \Gvf_2(\rho_0) - \Gvf_2(0), \Gvf_3(\rho_0)-\Gvf_3(0) ].
$$
Note that $b=0$ and $A$ is anisotropic.

The following problem arises naturally:

\begin{itemize}
\item[]
{\bf An over-determined problem for confocal ellipsoids}. {\sl Prove that if (ODP) admits a solution (in $H^1(\GO \setminus D)$), then $\GO$ and $D$ are confocal ellipsoids (or ellipses) and the common foci (when the volumes are fixed) is determined by the eigenvalues of $A$.}
\end{itemize}

The two-dimensional problem can be solved using the conformal mapping between $\GO \setminus \ol{D}$ and an annulus \cite[Theorem 10, p. 255]{Alfors}. If $A$ is isotropic, this problem is solved as the following theorem shows. The case of anisotropic $A$ has not been solved and is open.

\begin{thm}[\cite{KLS}]\label{conBall}
Let $D$ and $\GO$ be bounded  domains with Lipschitz boundaries in $\Rbb^3$ with $\overline{D} \subset \GO$. Suppose that $\GO \setminus \overline{D}$ is connected.
If (ODP) admits a solution for $A=c I$ for some constant $c$ where $I$ is the identity matrix in three dimensions, then $D$ and $\GO$ are concentric balls.
\end{thm}

As an immediate consequence, we have the following Corollary:

\begin{cor}
Suppose that $\Gs_c>\Gs_s$ and $\Gs_m$ is isotropic in addition to hypotheses of Theorem \ref{conBall}. If $(D, \GO)$ is neutral to multiple uniform fields, then $D$ and $\GO$ are concentric balls in three dimensions.
\end{cor}

Theorem \ref{conBall} is proved as follows.  Suppose that (ODP) admits a solution $w$ for $A = cI$.  Then, by (ODP), $|\Omega \setminus D| = -3c|D|$, and hence $c \not=0$. Introduce the angular derivative:
$$
A_{ij} = (x_j  + \frac{b_j}{c}) \p_i - (x_i  + \frac{b_i}{c}) \p_j, \quad i\not= j,
$$
where $b=(b_1,b_2,b_3)$ is the constant vector appearing in (ODP) and $\p_j$ denotes the partial derivative with respect to $x_j$-variable. Then one can see that $\Delta A_{ij}w = A_{ij}\Delta w = 0$ in $\Omega \setminus \overline{D}$, $A_{ij}w = 0$ on $\p\Omega$,  and $A_{ij} w = 0$ on $\p D$ provided that $A = cI$. Hence $A_{ij} w = 0$ in $\Omega\setminus\overline{D}$, which implies that $w$ is radial. Using this fact, one can prove Theorem \ref{conBall}. We emphasize that this argument using the angular derivative does not work if $A$ is not isotropic.

\subsection{The Newtonian potential formulation of the problem}

\begin{figure}[t!]
\begin{center}
\epsfig{figure=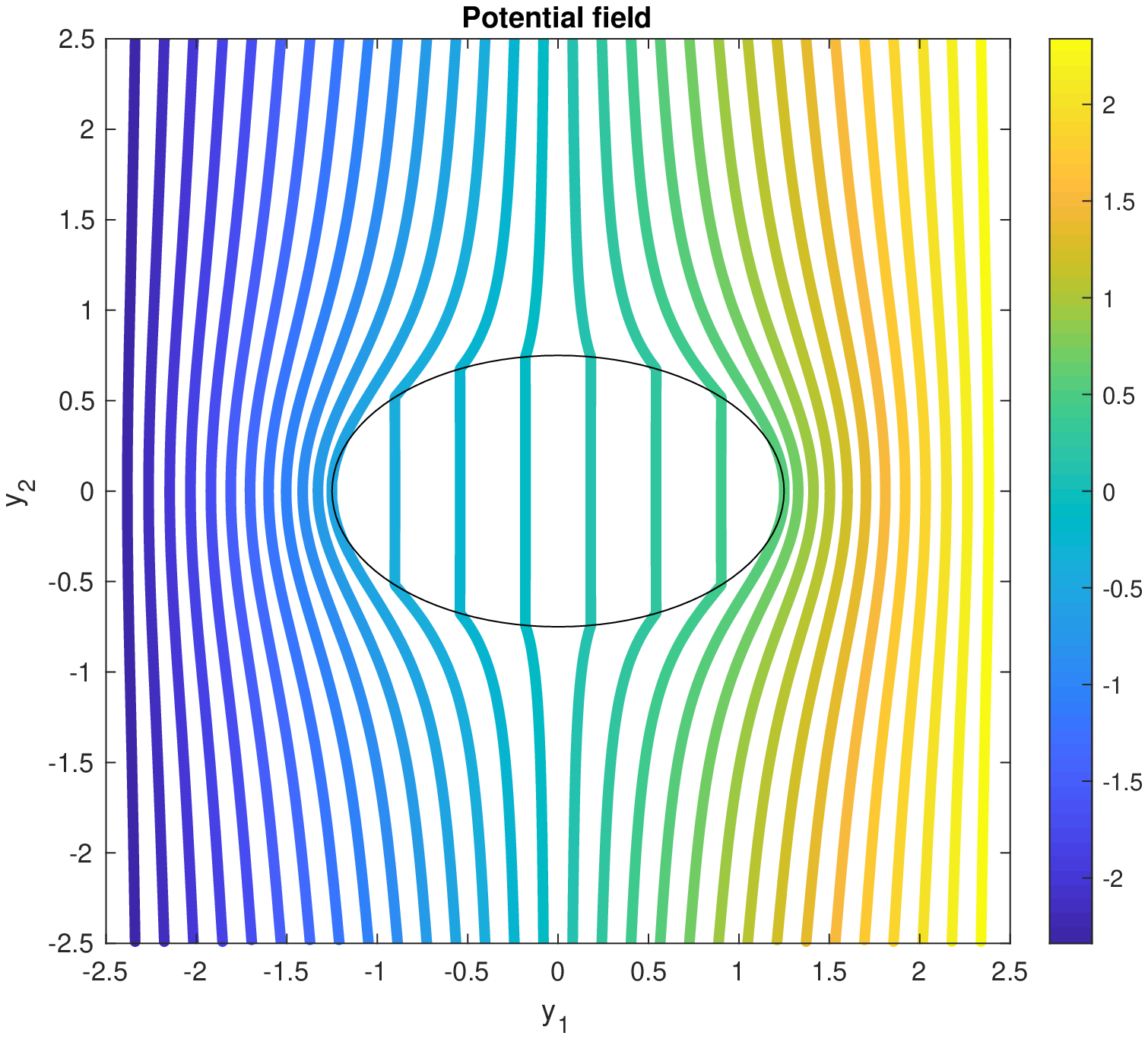, height=4cm} \epsfig{figure=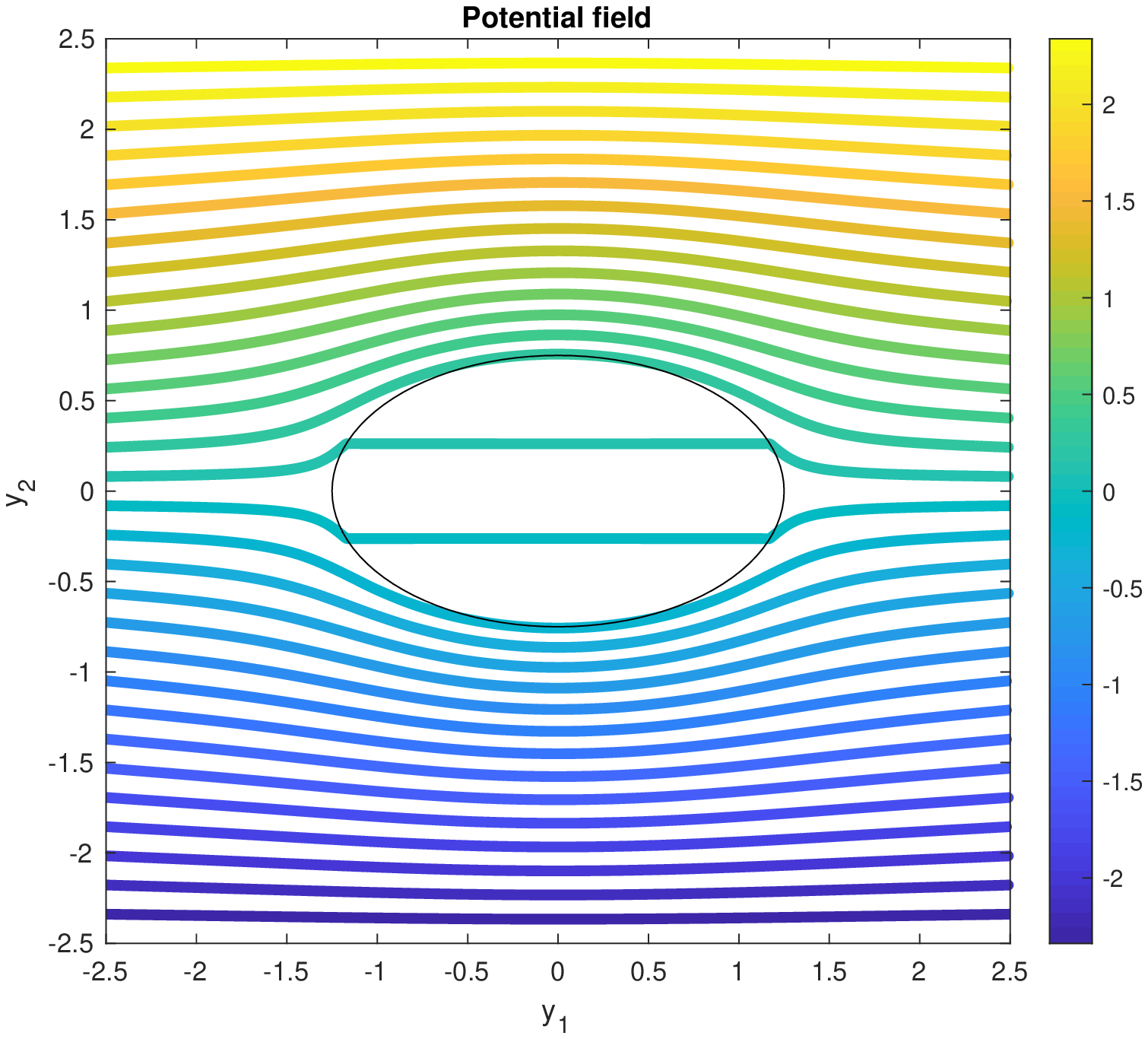, height=4cm}
\end{center}
\caption{Field inside an ellipse or an ellipsoid is uniform}\label{ellunif}
\end{figure}

Consider the conductivity problem (CP) when the conductivity distribution $\Gs$ is given by \eqnref{Gssimply}.
As one can see from Fig. \ref{ellunif}, the field inside $D$ is uniform if $D$ is an ellipse (or an ellipsoid). This a rather surprising fact that the field inside elliptic or
ellipsoidal inclusions is uniform seems to have been known for long time and
its proofs go back to Poisson (1826) and Maxwell (1873) (see \cite{Mac}). The converse is also true as we explain it in the sequel. For doing so, we need to recall the notion of the Newtonian potential.

The Newtonian potential of the domain $D$, which we denote by $N_D$, is defined by
\beq
N_D(x):= \frac{1}{|D|} \int_{D} \GG(x-y) dy .
\eeq
Usually the Newtonian potential is defined without the averaging factor $1/|D|$, but here it is more convenient to define it with the averaging factor. Since $\GD N_D(x) = 1/|D|$ for $x \in D$, we have
$$
N_D= \mbox{a quadratic part + a harmonic part} \quad\mbox{in } D.
$$
If $D$ is an ellipse or an ellipsoid, then the harmonic part of $N_D$ is quadratic and so is $N_D$ inside $D$. In fact, this is equivalent to the fact that the field inside elliptic or
ellipsoidal inclusions is uniform. Moreover, this property of the Newtonian potential's being a quadratic function inside the domain characterizes the ellipsoid and the ellipse:

\begin{thm}
Let $D$ be a simply connected bounded domain with the Lipschtiz boundary. If $N_D$ is quadratic inside $D$, then $D$ is an ellipse or an ellipsoid.
\end{thm}

This characterization of ellipsoids was proved by Dive in
1931 \cite{div31} and by Nikliborc in 1932 \cite{nik32} (see also
\cite{df86}). The reason why Dive and Nikliborc considered this
problem was to prove the converse of a theorem of Newton. Let $D$
be a simply connected domain whose center of mass is $0 \in D$ and
let $\Gl D$ be a dilation of $D$ by $\Gl >1$, {\it
i.e.}, $\Gl D =\{ \Gl x : x \in D \}$. A theorem due to Newton states that if $D$ is an
ellipsoid, then the gravitational force induced by the mass $\Gl D \setminus D$ is zero in $D$
\cite{kellog}. Dive and Nikliborc independently proved that the converse is true: If the
gravitational force induced by the uniform mass on $\Gl D
\setminus D$ is zero in $D$, then $D$ must be an
ellipsoid.

The following theorem was proved using the characterization of ellipsoids by Newtonian potentials.

\begin{thm}[\cite{KM06, KM08}]\label{thm:ell}
The following are equivalent for a simply connected bounded domain $D$:
\begin{itemize}
\item[{\rm (i)}] The polarization tensor $M$ associated with $D$ attains the lower Hashin-Shtrikman bound \eqnref{HSlower}.
\item[{\rm (ii)}] The solution to the conductivity problem (CP) when the conductivity distribution $\Gs$ is given by \eqnref{Gssimply} is linear inside $D$.
\item[{\rm (iii)}] $D$ is an ellipse in two dimensions and an ellipsoid in three dimensions.
\end{itemize}
\end{thm}

This theorem proves the generalized P\'olya-Szeg\"o conjecture explained before. That the linearity of the solution to (CP) when $\Gs$ is given by \eqnref{Gssimply} characterizes ellipsoids is known as the Eshelby's conjecture in the field of composites theory. Actually the Eshelby's conjecture (1961) \cite{esh61} asserts that the inclusion inside which the field is
uniform (or equivalently, the strain is constant) for a uniform loading is an ellipse or an ellipsoid. The corresponding  conjecture for the electro-static case is proved by Theorem \ref{thm:ell}. The Eshelby's conjecture (for the elasto-static case) was proved by Sendecyj \cite{sen70} in two dimensions and by Kang-Milton \cite{KM08} and Liu \cite{liu07} in three dimensions.

We now formulate the over-determined problem for the confocal ellipsoids in terms of the Newtonian potential.  It is proved in \cite{KLS} that the problem (ODP) admits a solution if and only if
\beq\label{11}
N_\GO(x)- N_D(x) =
\begin{cases}
0, \quad & x \in \Rbb^3 \setminus \GO, \\
\mbox{a quadratic polynomial}, \quad & x \in D.
\end{cases}
\eeq
Now the problem is to show that $D$ and $\GO$ are confocal ellipsoids if \eqnref{11} holds. If $\GO$ and $D$ are confocal ellipsoids, then both $N_\GO$ and $N_D$ are quadratic polynomials inside $D$, and so is $N_\GO- N_D$. A proof of the fact that $N_\GO= N_D$ outside $\GO$ can be found in \cite[p.61]{Mac}.
In the problem \eqnref{11}, the quadratic polynomial inside $D$ determines the common foci of $D$ and $\GO$. For example, one can see from Theorem \ref{conBall} that if the quadratic polynomial is of the form $c |x|^2+ \mbox{l.o.t}$, then $\GO$ and $D$ are concentric balls.

Now we can show the following theorem.
\begin{thm}\label{thm:equiv}
Suppose that $\p D$ is connected and $\Rbb^d \setminus \overline{D}$ is simply connected. Consider the following statements:
\begin{itemize}
\item[{\rm (i)}] $(D, \GO)$ is neutral to multiple uniform fields for some $\Gs$ given by \eqnref{Gsdouble}.
\item[{\rm (ii)}] The problem (ODP) admits a solution for some $A$ and $b$.
\item[{\rm (iii)}] The Newtonian potential formulation \eqnref{11} holds.
\end{itemize}
The following implications hold to be true:
\beq\label{threeimpli}
\mbox{{\rm (i)} $\Rightarrow$ {\rm (ii)} if $\Gs_c>\Gs_s$ }, \quad \mbox{{\rm (ii)} $\Rightarrow$ {\rm (iii)}}, \quad \mbox{{\rm (iii)} $\Rightarrow$ {\rm (i)}}.
\eeq
\end{thm}
\pf
The first implication was proved in \cite{KLS} and the proof is briefly reviewed at the beginning of this section. The second implication was proved in the same paper. We prove the third implication, namely, if \eqnref{11} holds, then there are conductivities $\Gs_c$, $\Gs_s$ and $\Gs_m$ such that $(D, \GO)$ is neutral to multiple uniform fields.

Let
\beq\label{defw}
w(x) := |\GO| N_\GO(x)- |\GO| N_D(x).
\eeq
By a rotation and a translation, if necessary, we may assume that
\beq \label{soltoODP}
w(x) = \begin{cases}
	0, \quad & x \in \Rbb^3 \setminus \GO, \\
	\sum_{j=1}^3 \Ga_j x_j^2 + \Ga , \quad & x \in D,
\end{cases}
\eeq
for some constants $\Ga_1, \Ga_2, \Ga_3$ and $\Ga$. In particular, there is no linear term in the quadratic function.
Define $u_j$ by
\beq
u_j(x) := \Gb_j \p_j \left( w(x) + \frac{1}{2} \sum_{j=1}^3 \Gb_j^{-1} x_j^2 \right),
\eeq
where $\Gb_j$'s are defined by \eqref{beta} with the conductivities to be determined later.

We claim that $u_j$ is the solution to (CP) with \eqref{neutral} when $a=e_j$. In fact, we see from the definition \eqnref{defw} of $w$ that
\begin{align*}
\p_j w(x) &= \int_\GO \p_{x_j} \GG(x-y) \;dy - f^{-1} \int_D \p_{x_j} \GG(x-y) \;dy \\
&= - \int_\GO \p_{y_j} \GG(x-y) \;dy + f^{-1} \int_D \p_{y_j} \GG(x-y) \;dy \\
&= - \int_{\p \GO} \GG(x-y) \;\Gv_j(y) \;dS(y) + f^{-1} \int_{\p D}  \GG(x-y) \; \Gv_j(y)\;dS(y),
\end{align*}
where the last equality follows from the divergence theorem. Here, $f$ is the volume fraction, namely, $f=|D|/|\GO|$. Thus,
\beq\label{pjw}
\p_j w(x)= - \Scal_{\p \GO} [\Gv_j](x) + f^{-1} \Scal_{\p D} [\Gv_j](x).
\eeq
Since the single layer potential is continuous across the boundary, $u_j$ is continuous across the interfaces $\p \GO$ and $\p D$.

Thanks to the jump relation \eqnref{jump} and \eqref{soltoODP}, we have, on $\p\GO$,
$$
\p_\nu (\p_j w)|_+ = -  \left( \frac{1}{2}I + \Kcal^*_{\p \GO} \right) [\Gv_j] + f^{-1} \p_\Gv \Scal_{\p D} [\Gv_j] =0,
$$
and hence
$$
\p_\nu (\p_j w)|_- = - \left( -\frac{1}{2}I + \Kcal^*_{\p \GO} \right) [\Gv_j] + f^{-1} \p_\Gv \Scal_{\p D} [\Gv_j] =\Gv_j.
$$
Thus,
$$
\Gs_{m,j} \p_\Gv u_j|_+ = \Gs_{m,j} \Gv_j,
$$
and
$$
\Gs_{s} \p_\Gv u_j|_- = \Gs_s (\Gb_j + 1 ) \Gv_j ,
$$
on  $\p \GO$. Since $\Gs_{m,j}= \Gs_s (\Gb_j + 1 )$ by the definition \eqref{beta} of $\Gb_j$,  we infer that
\beq
\Gs_{m,j} \p_\Gv u_j|_+  = \Gs_{s} \p_\Gv u_j|_- \quad\mbox{on } \p \GO.
\eeq

Similarly, thanks to \eqref{jump}, we have from \eqref{soltoODP} and \eqnref{pjw} that, on $\p D$,
$$
\p_\nu (\p_j w)|_- = -  \p_\Gv \Scal_{\p \GO} [\Gv_j] + f^{-1} \left( -\frac{1}{2}I + \Kcal^*_{\p D} \right)[\Gv_j] = 2\Ga_j \Gv_j,
$$
and hence
$$
\p_\nu (\p_j w)|_+ = - \p_\Gv \Scal_{\p \GO} [\Gv_j] + f^{-1} \left( \frac{1}{2}I + \Kcal^*_{\p D} \right)[\Gv_j] = (2\Ga_j+f^{-1}) \Gv_j.
$$
Thus, we have
$$
\Gs_{c} \p_\Gv u_j|_- = \Gs_c \Gb_j (2\Ga_j + \Gb_j^{-1}) \Gv_j
$$
and
$$
\Gs_{s} \p_\Gv u_j|_+ = \Gs_s \Gb_j (2\Ga_j + f^{-1} + \Gb_j^{-1}) \Gv_j.
$$
Thus,
\beq
\Gs_{s} \p_\Gv u_j|_+ = \Gs_{c} \p_\Gv u_j|_- \quad \text{ on } \p D
\eeq
if and only if
$$
\Gs_c \Gb_j (2\Ga_j + \Gb_j^{-1}) = \Gs_s \Gb_j (2\Ga_j + f^{-1} + \Gb_j^{-1}),
$$
or equivalently, by letting $\Gg:= 1-\Gs_c/\Gs_s$,
\beq\label{effective_general}
2\Ga_j \Gb_j \Gg + f^{-1} \Gb_j + \Gg = 0.
\eeq
So if we choose $\Gg$ and $\Gb_j$ (or $\Gs_c$, $\Gs_s$ and $\Gs_{m,j}$) so that \eqref{effective_general} holds, then $(D,\GO)$ is neutral to multiple uniform fields.

There is yet another restriction when we solve \eqref{effective_general} for $\Gg$ and $\Gb_j$; $\Gs_c$, $\Gs_s$ and $\Gs_{m,j}$ should be positive. This condition can be easily fulfilled. In fact, the following relation follows easily from \eqref{effective_general}:
$$
\Gs_{m,j} = \Gs_s \left( 1- \frac{\Gg}{2\Ga_j \Gg + f^{-1}} \right).
$$
The quantity $2\Ga_j \Gg + f^{-1}$ in the above is nonzero since $\Gg$ can be chosen small as we see shortly. Thus the positivity is achieved if
$$
1- \frac{\Gg}{2\Ga_j \Gg + f^{-1}} >0
$$
which in turn can be achieved by taking $\Gg$ so that
$$
|\Gg| \le \min_{1 \le j \le 3} \frac{f^{-1}}{|1-2\Ga_j|+1}.
$$
This completes the proof.
\qed

In the course of the proof, we derived the neutrality condition for confocal ellipsoids.
\begin{cor}
Let $D$ and $\GO$ be confocal ellipsoids whose boundaries are respectively given by \eqnref{ellD} and \eqnref{ellGO}.     If the conductivity distribution $\Gs$ given by \eqnref{Gsdouble} and \eqnref{Gsmaniso} satisfies
\beq\label{neutralell}
2\Ga_j f\left( \frac{\Gs_{m,j}}{\Gs_s} - 1\right)\left( 1 - \frac{\Gs_c}{\Gs_s} \right) + \left( \frac{\Gs_{m,j}}{\Gs_s} - 1\right) + f \left( 1 - \frac{\Gs_c}{\Gs_s} \right) = 0, \;\; j=1,2,3,
\eeq
then $(D, \GO)$ is neutral to multiple uniform fields. Here $f$ is the volume fraction and $\Ga_j$'s are constants given by \eqnref{soltoODP}, 	\textit{i.e.},
\beq
\Ga_j = -\frac{1}{4}\int_0^{\Gr} \frac{1}{c_j^2 + s} \frac{\sqrt{(c_1^2 + \Gr)(c_2^2 + \Gr)(c_3^2 + \Gr)}}{\sqrt{(c_1^2 + s)(c_2^2 + s)(c_3^2 + s)}} \;\; ds, \; j=1,2,3.
\eeq
\end{cor}

Let us look into the neutrality condition \eqref{effective_general} or \eqnref{neutralell} further. According to \eqref{defw} and \eqref{soltoODP},
$$
1-f^{-1} = \GD w = 2\sum_{j=1}^3 \Ga_j \quad\mbox{in } D.
$$
We then have from \eqref{effective_general}
$$
1-f^{-1} = - \frac{3f^{-1}}{\Gg} - \sum_{j=1}^3 \frac{1}{\Gb_j} ,
$$
and hence
$$
-1 + \frac{3}{\Gg} + f \left( 1+ \sum_{j=1}^3 \frac{1}{\Gb_j} \right)=0.
$$
Writing it in terms of conductivities, we have
\beq
\frac{2\Gs_s+\Gs_c}{\Gs_s-\Gs_c} + \frac{f}{3} \sum_{j=1}^3 \frac{\Gs_{m,j} + 2\Gs_s}{\Gs_{m,j} -\Gs_s} =0.
\eeq
This is a necessary neutrality condition when $\Gs_m$ is anisotropic. In particular, if $\Gs_m$ is a scalar, namely, $\Gs_{m,j}=\Gs_m$, then it is exactly the neutrality condition \eqref{effective} of concentric balls.

\subsection{Quadrature domains-Neumann ovaloids}\label{subsec:ovaloid}

Let us look further into the Newtonian potential formulation \eqref{11} of the problem. The problem is to prove that if it holds, then $D$ and $\GO$ must be confocal ellipsoids. We show that the condition \eqref{11} in $\Rbb^3 \setminus \GO$ alone does not yield the answer.

The condition \eqref{11} in $\Rbb^3 \setminus \GO$ yields
$$
\int_{\p\GO} {N}_{\GO}(x) g(x) \, dS = \int_{\p\GO} {N}_{D}(x) g(x) \, dS
$$
for any $g \in H^{-1/2}(\p\GO)$. By changing the order of integrations, we have
\beq\label{mean}
\frac{1}{|\GO|} \int_{\GO} u(x) dx = \frac{1}{|D|} \int_{D} u(x) dx,
\eeq
where $u(x)= \Scal_{\p\GO}[g](x)$, $x\in \GO$. Thus \eqnref{mean} holds for all $H_h^1(\GO)$ where subscript $h$ means that it is a collection of harmonic functions in $\GO$.  The condition \eqnref{mean} does not guarantee that $D$ and $\GO$ are confocal ellipsoids as will be seen in what follows, and the condition \eqref{11} in $D$ should be utilized.

In fact, an open set $\GO \subset \Rbb^d$ is called a \emph{quadrature domain} (see, e.g., \cite[(4.1)]{Shapiro}, and also \cite{GS,Sakai}) if there exists a distribution $\Gm$ with a compact support in $\GO$ such that
\beq
\int_{\GO} u(x) dx = \la \mu, u\ra \quad \text{ for all } u \in H_h^1(\GO).
\eeq
The simplest quadrature domain may be balls: It is well known as the mean value property:
\beq
\int_{\GO} u(x)  dx = |\GO| u(c) \quad \text{ for all } u \in H_h^1(\GO),
\eeq
where $c$ is the center of the ball. In this case the distribution $\Gm$ is the point mass (the Dirac delta) multiplied by the volume of $\GO$.

Ellipsoids are also quadrature domains. Let
$$
\GO= \left\{ x \in \Rbb^{d} \; | \; \sum_{i=1}^{d} \frac{x_i^2}{a_i^2}  < 1 \right\} \quad (a_1 \geq a_2 \geq \cdots \geq a_{d-1} > a_d >0),
$$
and let
$$
F= \left\{ x \in \Rbb^{d-1} | \; \sum_{i=1}^{d-1} \frac{x_i^2}{a_i^2 - a_d^2} <1 \right\}.
$$
The lower dimensional set $F$ is called the focal ellipsoid of $\GO$.  The following quadrature identity holds (see, e.g., \cite{Khavinsion, KhaL}):
\beq\label{MVPell}
\int_{\GO} u(x) dx = 2\prod_{i=1}^{d-1}\frac{a_i}{(a_i^2 - a_d^2)^{1/2}} \int_F \left(1 - \sum_{i=1}^{d-1} \frac{x_i^2}{a_i^2 - a_d^2} \right)^{1/2} \; u(x',0) \;dx'
\eeq
for all $u \in H_h^1(\GO)$. Here, $x'$ is $(x_1, \dots, x_{d-1})$. Note that if $D$ and $\GO$ are confocal ellipsoids, then their focal ellipsoids are the same, and hence \eqnref{mean} holds.

\begin{figure}[t!]
\begin{center}
\epsfig{figure=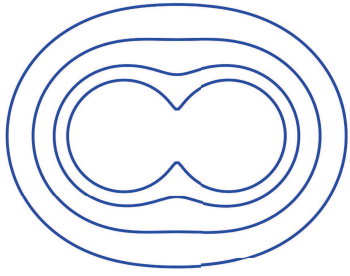, height=4cm}
\end{center}
\caption{Neumann ovals with same foci} \label{Noval}
\end{figure}

There is yet another class of domains satisfying \eqnref{mean}.
A domain $\GO \subset \Rbb^d$ is called a \emph{Neumann ovaloid} if it admits the following quadrature identity
\beq\label{NeumannOvaloid}
\int_{\GO} u(x) dx = C\left(u(p_1) + u(p_2) \right) \quad \text{ for all }  u \in H_h^1(\GO),
\eeq
where $p_1$ and $p_2$ are distinct points in $\Rbb^d$ and $C>0$ is a constant. If $C$ is sufficiently small compared to $|p_1 - p_2|$, then a union of two balls of the same radius centered at $p_1$ and $p_2$ satisfies the identity \eqref{NeumannOvaloid}. However, if $C$ is sufficiently large, then there is an axially-symmetric domain satisfying \eqref{NeumannOvaloid}. For example, if $\GO$ is the domain in $\Rbb^2$ bounded by the curve
\beq
\left( x_1^2 + x_2^2 \right)^2 = \Ga^2 \left( x_1^2 + x_2^2 \right) + 4\Gve^2 x_1^2,
\eeq
where $\Ga$ and $\Gve$ are some positive constants (see Fig. \ref{Noval}), then it admits a quadrature identity \eqnref{NeumannOvaloid} with $C= |\GO|/2$ (see \cite[p. 19--20]{Shapiro} for a proof). In this case, the relation among $|\GO|$, $\Ga$ and $\Ge$ is given by $|\GO| = \Gp (\Ga^2 + 2\Gve^2)$.
These two-dimensional Neumann ovals were discovered by C. Neumann \cite{Neumann}. The existence and uniqueness of the higher dimensional Neumann ovaloids are known (see \cite{KarpL} and references therein), but there is no known explicit expression except four-dimensional (and two-dimensional) ones, to the best of our knowledge. We refer to a recent paper \cite{KarpL} for an explicit parametrization of a four-dimensional Neumann ovaloid. If $(D,\GO)$ are Neumann ovaloids with same foci, then \eqnref{mean} holds.

\section{Weakly neutral inclusions}\label{sec:weak}

We now consider the weakly neutral inclusion problem presented at the end of Introduction, namely, the problem of coating a given domain of general shape by another domain so that the resulting inclusion of core-shell structure satisfies the weak neutrality condition \eqnref{wneutral}, namely, its polarization tensor vanishes. In this section we present two classes of domains for which the weakly neutral inclusion problem can be solved. The first one is defined by a conformal mapping from the exterior of the unit disk onto the exterior of the domain, and construction of the coating is explicit. This result is new. The other class are small perturbations of a disk for which existence of a coating is proved. This result is from \cite{KLS2D, KLS3D}. Note that neutral inclusions are weakly neutral to multiple uniform fields. Thus concentric disks or balls can be realized as weakly neutral inclusions. However, no other examples of weakly neutral inclusions are known.

\subsection{$b_D$-vanishing domains}

Let $D$ be a bounded domain in $\Cbb=\Rbb^2$ with the Lipschitz continuous boundary, and let $z=\Phi(\Gz)$ is the Riemann mapping from $|\Gz|>1$ ($\Cbb\setminus\ol{U}$, where $U$ is the unit disk) onto $\Cbb\setminus\ol{D}$.  The conformal mapping $\Phi$ takes the form
$$
\Phi(\Gz)=b_{-1}\Gz+b_0+\frac{b_1}{\Gz}+ \mbox{h.o.t}.
$$
By dilating and translating $D$ if necessary, we assume that $b_{-1}=1$ and $b_0=0$, and denote $b_1$ by $b_D$, so that the Riemann mapping $\Phi$ takes the form
\beq\label{f}
\Phi(\Gz)=\Gz+\frac{b_D}{\Gz}+\mbox{h.o.t}.
\eeq
The domains we consider in this subsection are such that $b_D=0$.

We want to coat $D$ by an another bounded domain $\GO$ so that the PT of the coated structure vanishes. Let the conductivity distribution $\Gs$ be given by \eqnref{Gsdouble}. Furthermore, we assume that $\Gs_c=0$ or $\infty$. This assumption is required because we use the conformal mapping from $\Cbb\setminus\ol{\GD}$ onto $\Cbb\setminus\ol{D}$. We assume that $\Gs_c=\infty$, and the other case can be handled in the same. We also assume $\Gs_m=1$ without loss of generality.

Since $\Gs_c=\infty$, (CP) is now of the form
$$
\mbox{(CP)}_\infty \quad
\begin{cases}
	\nabla \cdot \Gs \nabla u=0 \quad &\mbox{in}~\Rbb^2\setminus\ol{D},\\
	u=\Gl(\mbox{constant})  \quad&\mbox{on}~\p D,\\
	u(x)-a\cdot x=O(|x|^{-1}) \quad&\mbox{as}~|x|\rightarrow\infty,
\end{cases}
$$
where $\Gs= \Gs_s \chi(\GO \setminus D) + \chi(\Rbb^2 \setminus \GO)$. The constant $\Gl$ is determined by the condition $\int_{\p D}\ \p_\nu u|_{+}=0$. The problem is to find $\Gs_s$ and $\GO$ so that the solution $u$ to (CP)$_\infty$ satisfies the weak neutrality condition \eqnref{wneutral}.

Since $u(x)-a\cdot x$ tends to $0$ as $|x| \to \infty$ and $D$, $\GO$ are simply connected, there are functions $U_m$ and $U_s$ analytic in $\Cbb\setminus\ol{\GO}$ and $\GO\setminus\ol{D}$, respectively, such that $\Re U_m=u$ and $\Re U_s=u$ in their respective domains ($\Re$ stands for the real part). One can see using the Cauchy-Riemann equations that the transmission conditions to be satisfied by $u$ on $\p \GO$ is equivalent to
\beq\label{tr2}
(1+\Gs_s)U_s+(1-\Gs_s)\overline{U_s}=2U_m+c \quad\mbox{on } \p \GO,
\eeq
for some constant $c$. Moreover, $U_m$ admits the following expansion at $\infty$:
$$
U_m(z)=\Ga z+\frac{\Ga_1(\Ga)}{z}+\mbox{h.o.t},
$$
where $\Ga=a_1-ia_2$. Thus the weak neutrality condition \eqref{wneutral} is equivalent to
\beq\label{a1}
\Ga_1(\Ga)=0 \quad \mbox{for  all }  \Ga~(\mbox{or equivalently, for } \Ga=1,i).
\eeq

With the conformal mapping $\Phi$ in \eqnref{f}, let $V^\Ga=U\circ\Phi$. Then we have
\beq\label{Vzeta}
V^\Ga(\Gz)=\Ga\Phi(\Gz)+\frac{\Ga_1(\Ga)}{\Phi(\Gz)}+\mbox{h.o.t} =\Ga\Gz+\frac{\Ga b_D+\Ga_1(\Ga)}{\Gz}+\mbox{h.o.t}.
\eeq
Let $U'$ be a simply connected domain containing $\ol{U}$ defined by
\beq\label{U'def}
\Phi(\p U')= \p\GO.
\eeq
The transmission condition \eqref{tr2} is transformed by $\Phi$ to
\beq
(1+\Gs_s)V^\Ga_s+(1-\Gs_s)\overline{V^\Ga_s}=2V^\Ga_m+c \quad\mbox{on } \p U' .
\eeq

If $b_D=0$, then \eqnref{Vzeta} takes the form
$$
V^\Ga(\Gz)=\Ga\Gz+\frac{\Ga_1(\Ga)}{\Gz}+\mbox{h.o.t}.
$$
Thus \eqnref{a1} is fulfilled if and only if $V^\Ga$ satisfies
\beq\label{V2}
|V^\Ga(\Gz)-\Ga\Gz| = O\left(|\Gz|^{-2} \right) \quad\mbox{as } |\Gz|\rightarrow \infty,
\eeq
for $\Ga=1,i$. Since $\Re (V^\Ga)$ is the solution to (CP)$_\infty$ with $D$ and $\GO$ replaced by $U$ and $U'$, respectively, \eqnref{V2} is satisfied if $(U, U')$ is a neutral inclusion. Since $\Gs_c=\infty$ and $U$ is the unit disk, it suffices to take $U'$ to be a disk of radius $r$ and the conductivity $\Gs_s$ to satisfy the neutrality condition \eqnref{effective}, which is in this case
\beq\label{effective2}
(1 - \Gs_s) - r^2 (1 + \Gs_s) =0.
\eeq

In summary, for a domain $D$ such that $b_D=0$, we take $\Gs$ and $r$ to satisfy \eqnref{effective2}. Then $(D, \GO)$ where $\GO$ is define by \eqnref{U'def} is weakly neutral to multiple uniform fields.

\begin{figure}[h]
	\centering
	\begin{subfigure}{0.45\textwidth}
		\includegraphics[scale=0.3]{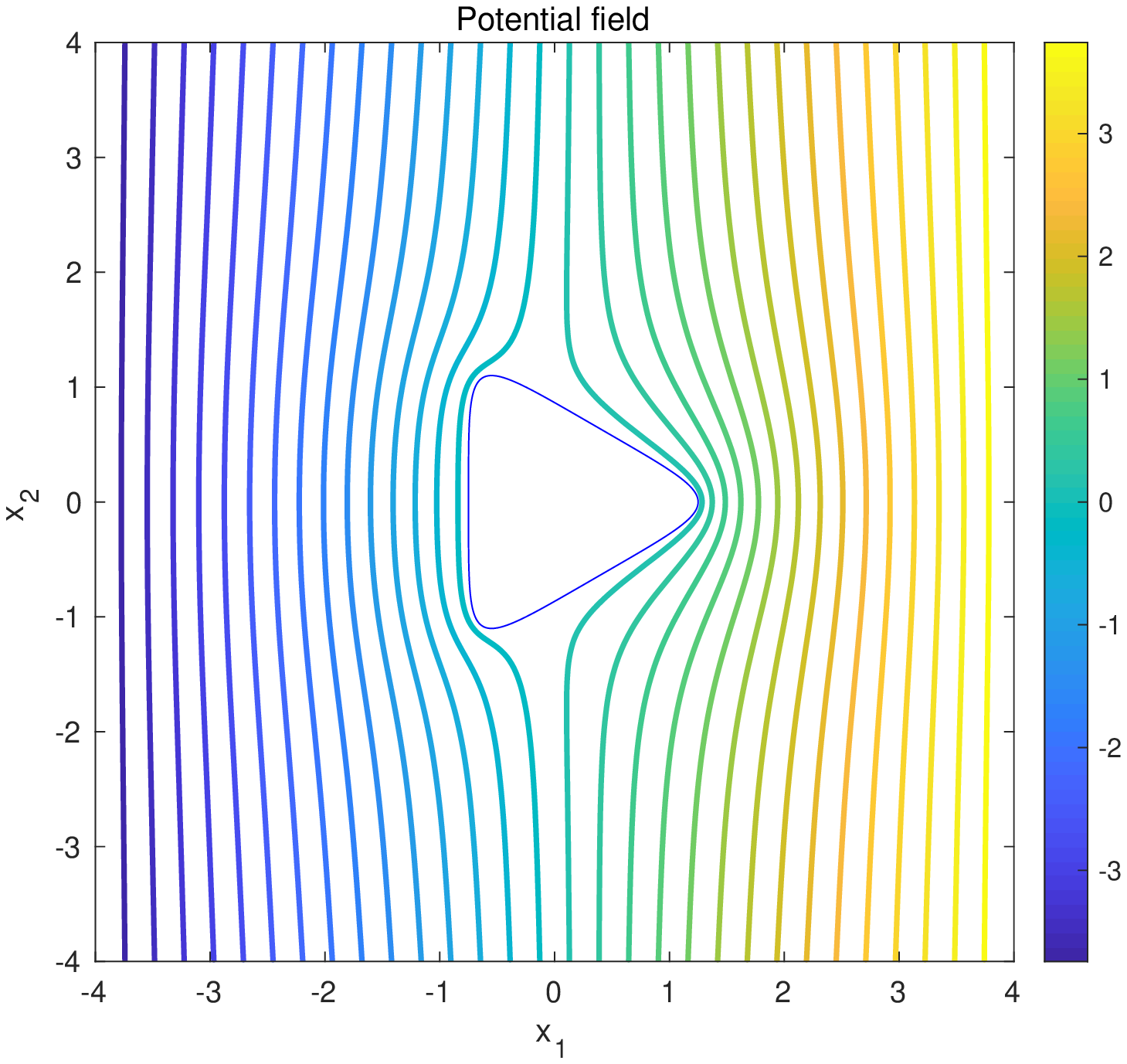}
		\subcaption{}
	\end{subfigure}
	\begin{subfigure}{0.45\textwidth}
		\includegraphics[scale=0.3]{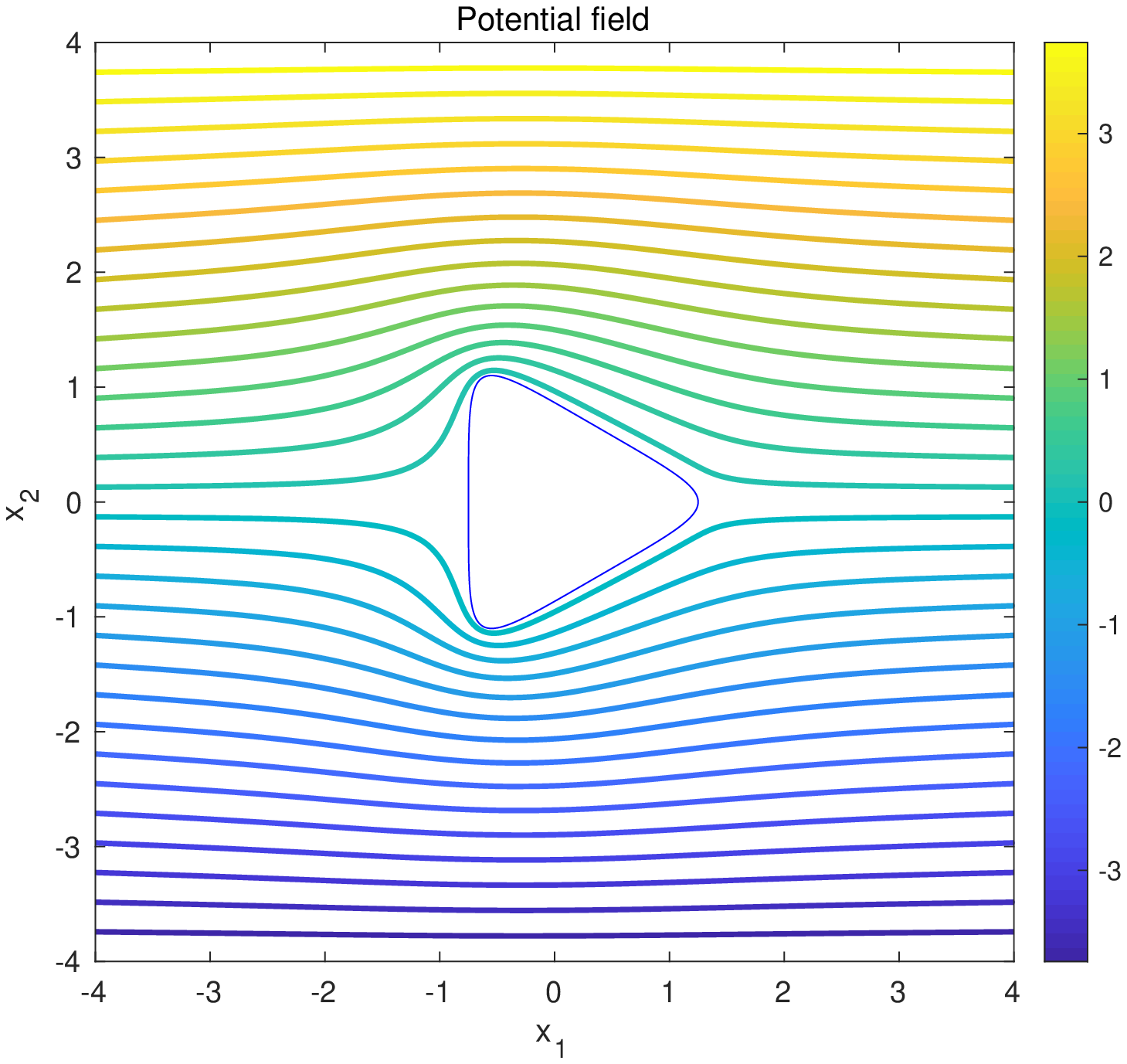}
		\caption{}
	\end{subfigure}\\
	\begin{subfigure}{0.45\textwidth}
		\includegraphics[scale=0.3]{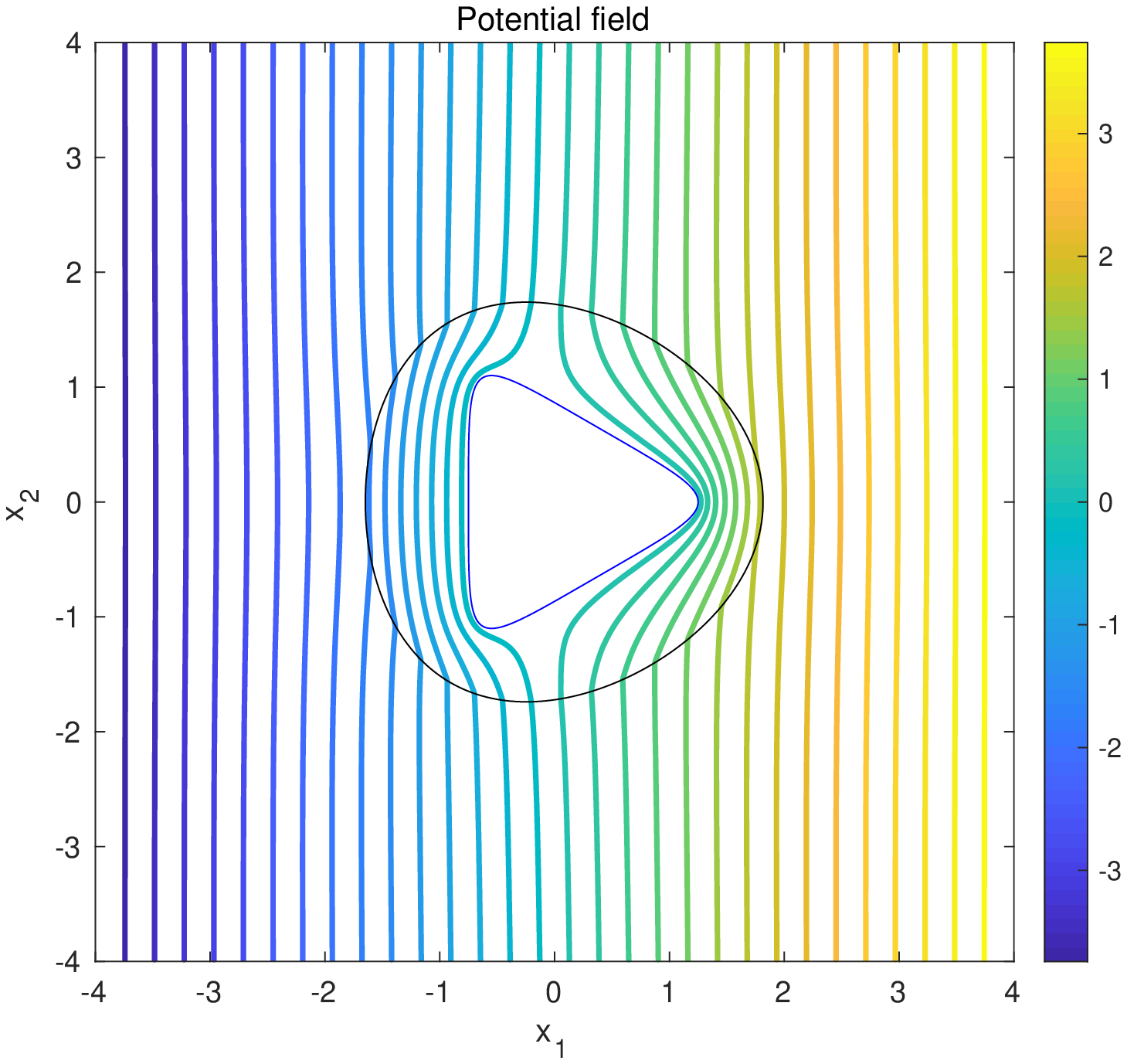}
		\caption{}
	\end{subfigure}
	\begin{subfigure}{0.45\textwidth}
		\includegraphics[scale=0.3]{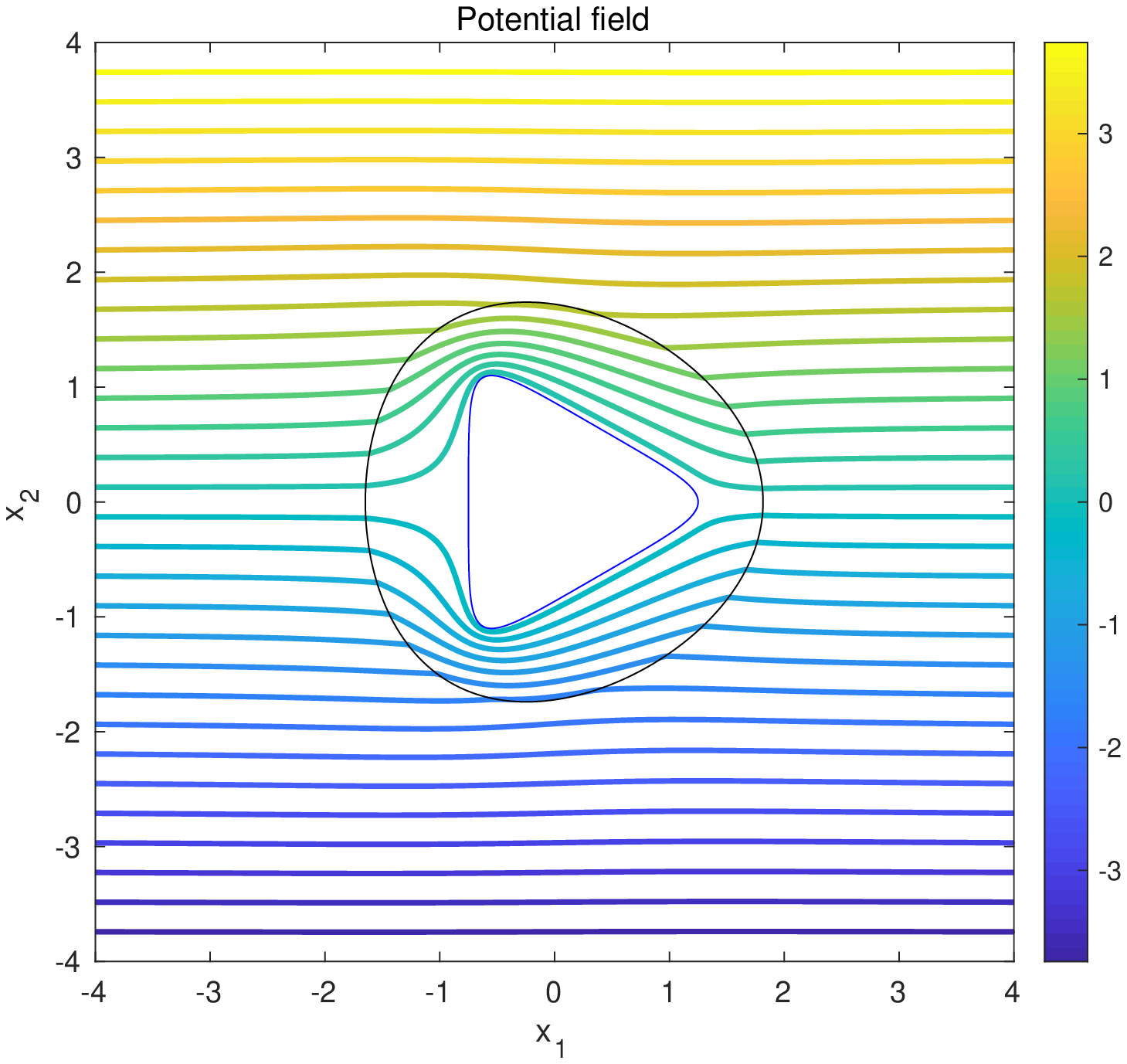}
		\caption{}
	\end{subfigure}
	\caption{The core-shell structure defined by the conformal mapping $\Phi(\Gz)=\Gz+\frac{1}{4\Gz^2}$. Field perturbation with the coating ((c) and (d)) is much weaker that with it ((a) and (b)). }
	\label{b2}
\end{figure}

\begin{figure}[h]
	\centering
	\begin{subfigure}{0.45\textwidth}
		\includegraphics[scale=0.3]{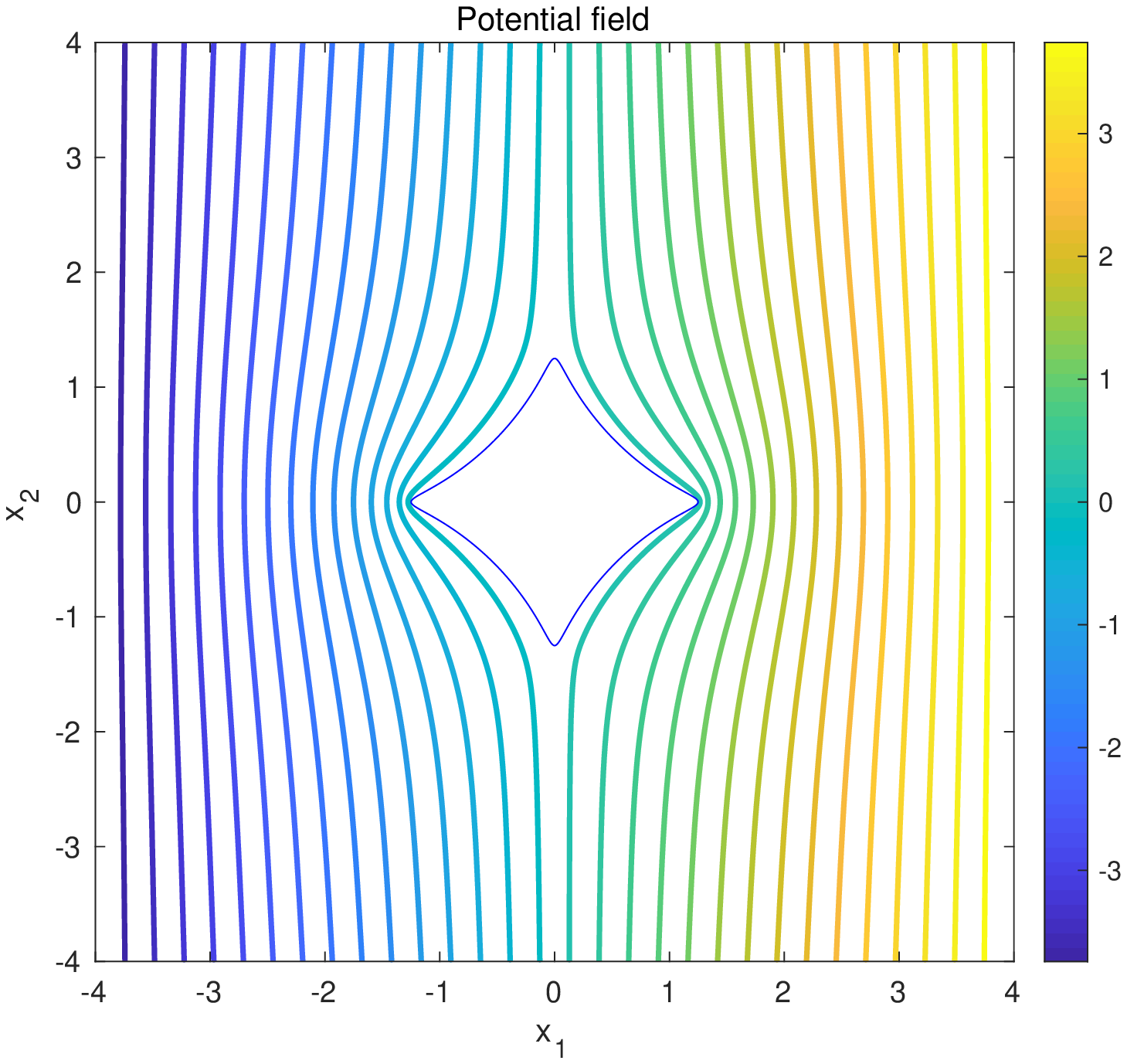}
		\subcaption{}
	\end{subfigure}
	\begin{subfigure}{0.45\textwidth}
		\includegraphics[scale=0.3]{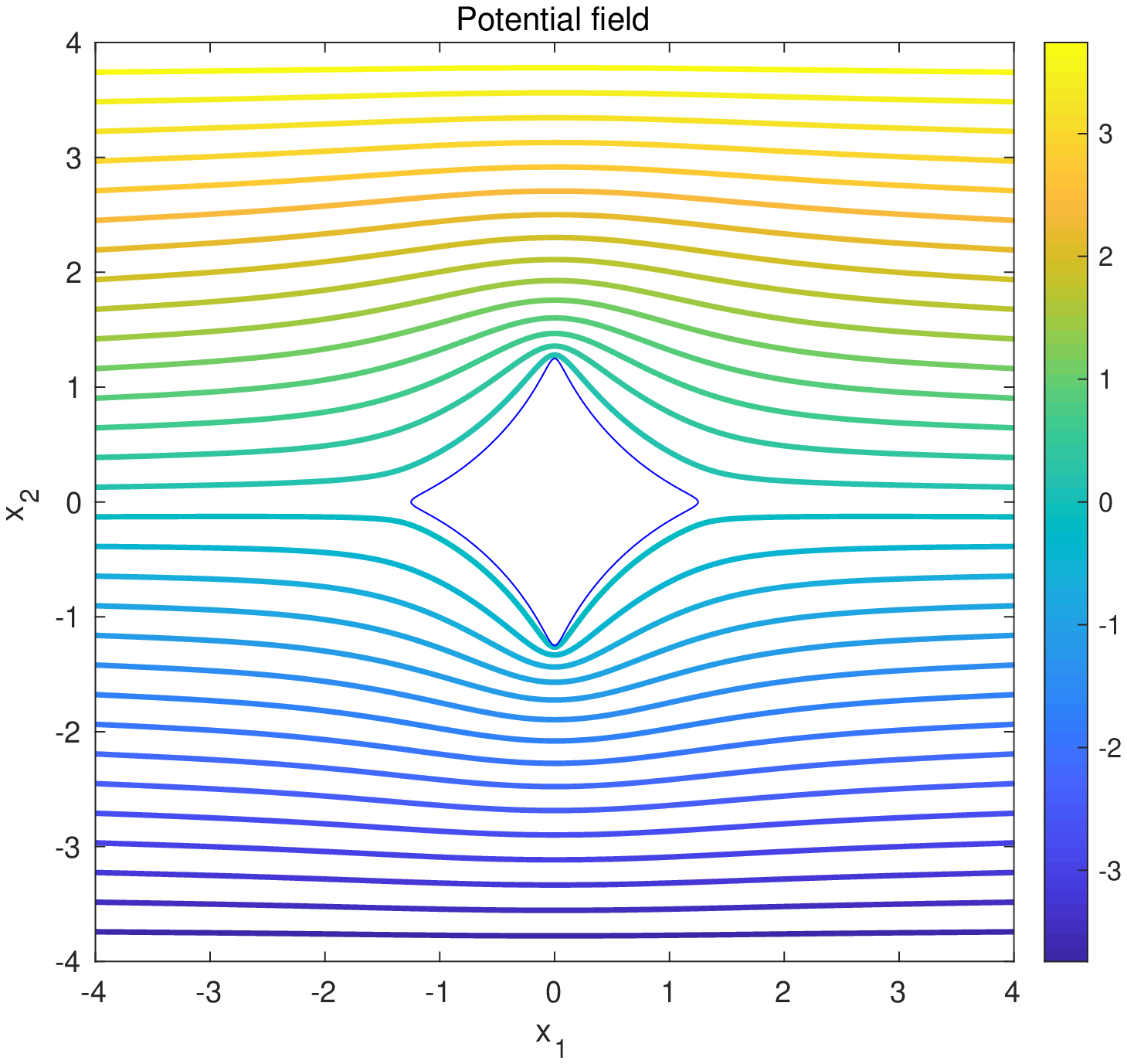}
		\caption{}
	\end{subfigure}\\
	\begin{subfigure}{0.45\textwidth}
		\includegraphics[scale=0.3]{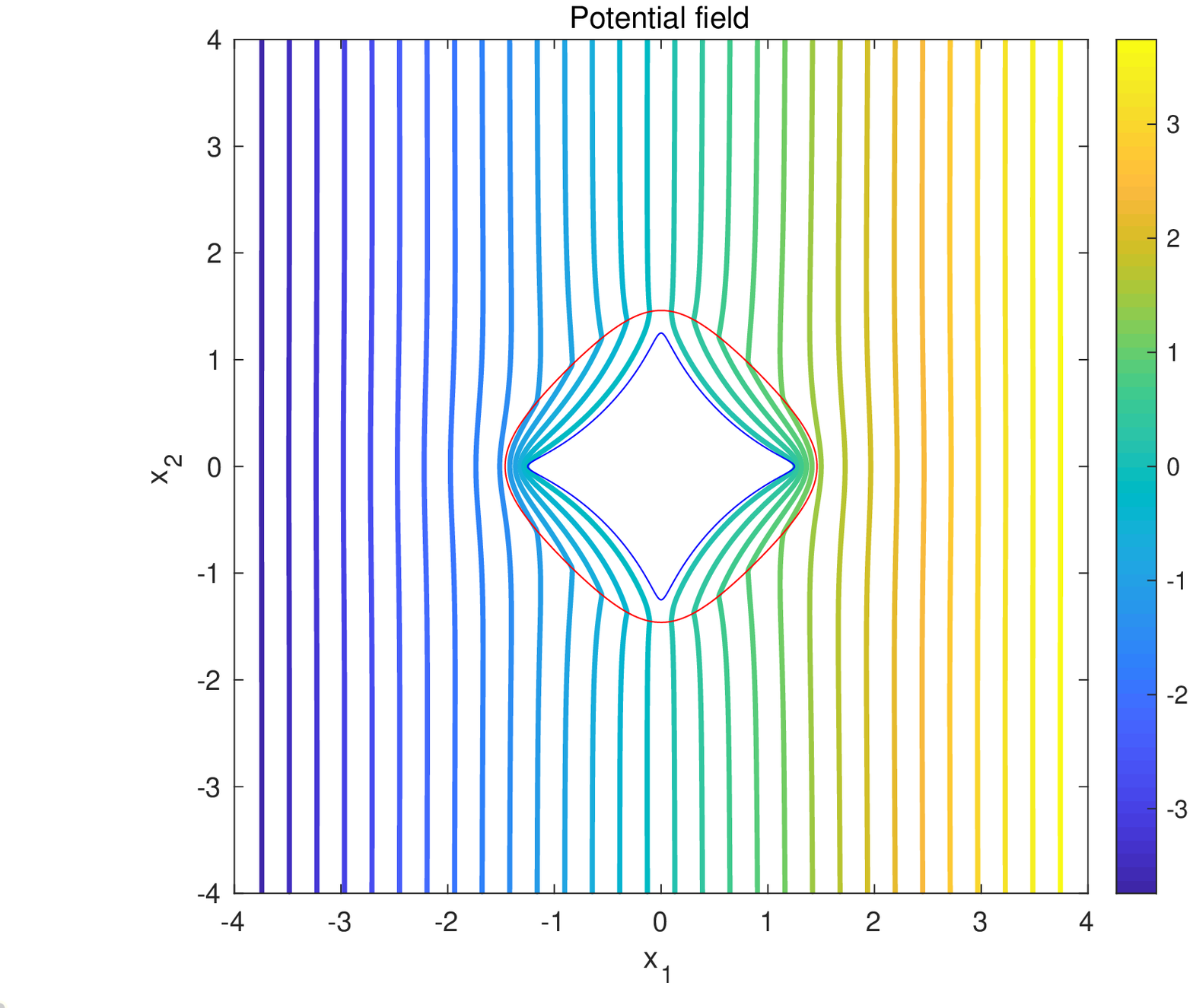}
		\caption{}
	\end{subfigure}
	\begin{subfigure}{0.45\textwidth}
		\includegraphics[scale=0.3]{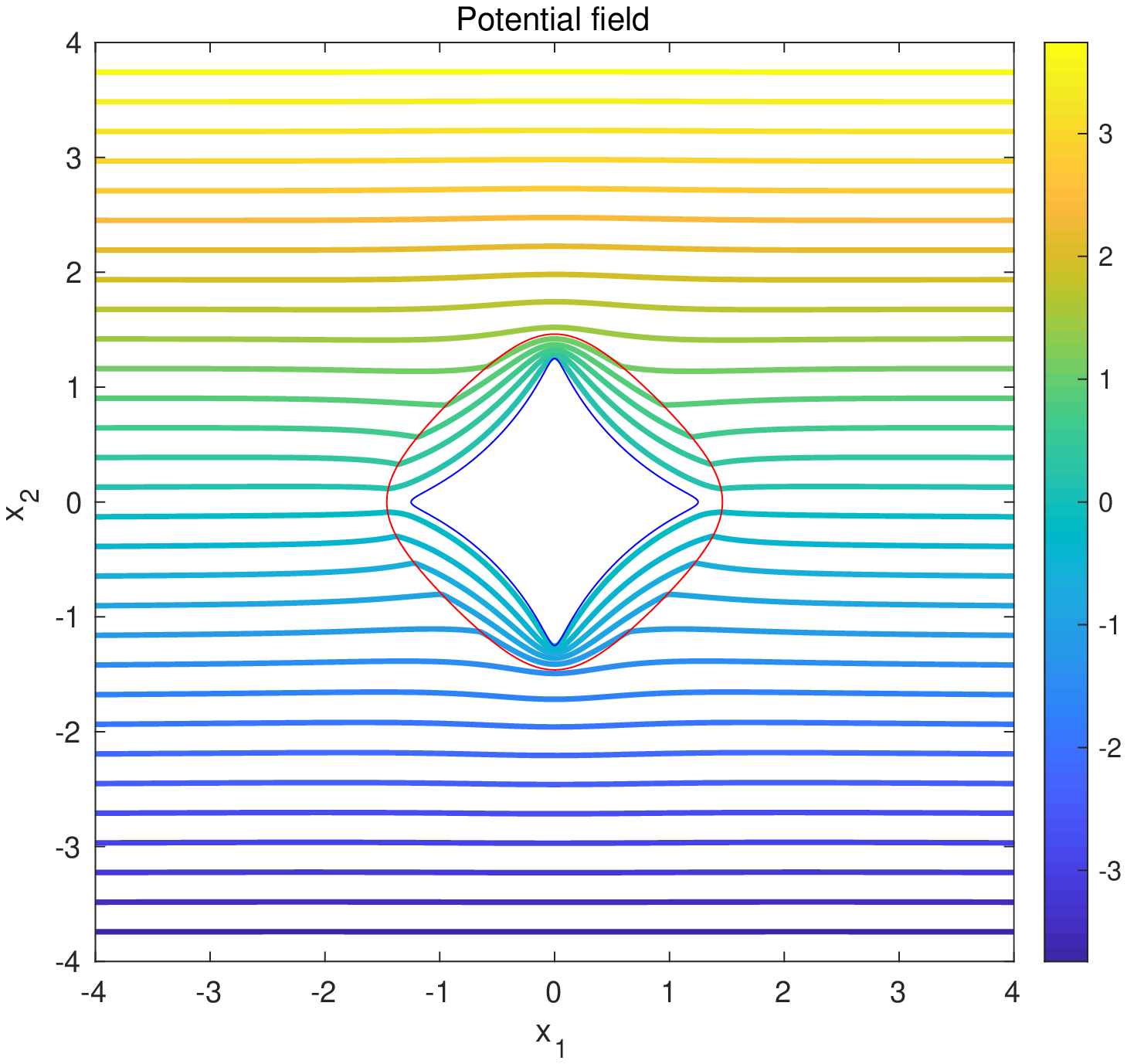}
		\caption{}
	\end{subfigure}
	\caption{The core-shell structure defined by the conformal mapping $\Phi(\Gz)=\Gz+\frac{1}{4\Gz^3}$. Field perturbation with the coating ((c) and (d)) is much weaker that with it ((a) and (b)). }
	\label{b3}
\end{figure}

We now present two results of numerical experiment. In Fig. \ref{b2}, the conformal mapping for the domain $D$ and the conductivity $\Gs_s$ are given by
\beq
\Phi(\Gz)=\Gz+\frac{1}{4\Gz^2} \quad\mbox{and}\quad \Gs_s=0.5.
\eeq
It shows the domains $D$ and its coating $\GO$ determined by the method described above. It clearly shows that field perturbation with the coating is much weaker than that without it. Fig. \ref{b3} is with the conformal mapping
\beq
\Phi(\Gz)=\Gz+\frac{1}{4\Gz^3} \quad\mbox{and}\quad \Gs_s=0.3.
\eeq

\subsection{Small perturbation of disks}

We now review the result from \cite{KLS2D} which shows that a small perturbation of a disk allows a coating such that the resulting core-shell structure is a weakly neutral inclusion, namely, a PT-cancelling structure. It is an existence result based on the implicit function theorem, so we do not know how small it can be.

Let $D_0$ be a disk of radius $r_i$ centered at the origin. For a given function $h$ on the unit circle $T$, the perturbation $D_h$ of $D_0$ is defined to be
\beq
\p D_h: =\left\{ ~x ~|~x =(r_i+ h(\hat{x}))\hat{x}, \quad |\hat{x}|=1 ~ \right\}.
\eeq
We consider $W^{2,\infty}(T)$ (the derivatives up to 2 are bounded) for a class of perturbation functions $h$.

To define domains for coating, we let $\GO_0$ be the disk of radius $r_e$ centered at the origin such that $(D_0, \GO_0)$ be a neutral inclusions, namely, the radius and the conductivities are chosen so that the neutrality condition \eqnref{effective} is satisfied. We then define perturbations of $\GO_0$ as follows:
\beq
\p \GO_b :=\left\{ ~x ~|~x =(r_e + b(\hat{x})) \hat{x}, \quad |\hat{x}|=1 ~ \right\},
\eeq
where $b$ is of the form
\beq\label{b}
b(\Gt)=b(\hat{x})=b_0+b_1\cos 2\Gt + b_2\sin 2\Gt.
\eeq
Here $b_0,b_1,b_2$ are real constants.

If $h$ and $b$ are sufficiently small, then $(D_h, \GO_b)$ defines an inclusion of the core-shell structure. Let $M(h,b)=M(D_h,\GO_b)$ be the PT of $(D_h, \GO_b)$ as defined in \eqnref{PT0}. Since $M$ is symmetric, we may regard $M$ as a three-dimensional vector-valued function. Since the collection of all $b$ of the form \eqnref{b} is of three dimensions, $M$ can be regarded as a mapping from $U \times V$ into $\Rbb^3$, where $U$ and $V$ are some neighborhoods of $0$ in $W^{2,\infty}(T)$ and $\Rbb^3$, respectively. Since $(D_0, \GO_0)$ is neutral, we have $M(0,0)=0$. It is then proved that
\beq\label{PTdet}
\mbox{det} \frac{\p M}{\p (b_0,b_1, b_2)} (0,0) \neq 0.
\eeq
Then an implicit function theorem is invoked to arrive at the following theorem.

\begin{thm}[\cite{KLS2D}]\label{2D}
There is $\Ge>0$ such that for each $h \in W^{2,\infty}(T)$ with $\| h \|_{2,\infty}<\Ge$ there is $b=b(h) \in \Rbb^3$ such that
\beq
M(h, b(h))=0,
\eeq
namely, the inclusion $(D_h,\GO_{b(h)})$ of the core-shell structure is weakly neutral to multiple uniform fields. The mapping $h \mapsto b(h)$ is continuous.
\end{thm}

Proving \eqnref{PTdet} is quite technical. This two-dimensional theorem has been extended in \cite{KLS3D} to three dimensions, which is even more technically complicated, to show that small perturbations of a sphere allow coatings so that the resulting inclusions of the core-shell structure are weakly neutral to multiple uniform fields. For that, the functions $b$ in \eqnref{b} is replaced with
\beq\label{b3D}
b(\hat{x})=b_0+ \sum_{j=1}^5 b_j Y_j^2 (\hat{x}),
\eeq
where $Y_j^2 (\hat{x})$ are spherical harmonics of order 2 (there are five linearly independent ones). Then the PT is regarded as a local mapping from $W^{2,\infty}(S) \times \Rbb^6$ ($S$ is the unit sphere) into $\Rbb^6$, and an analogy to \eqnref{PTdet} is proved.

\section{Weakly neutral inclusions by imperfect bonding}\label{sec:imperfect}

So far we consider neutral or weakly neutral inclusions of the core-shell structure. There is yet another method to achieve neutrality: It is by introducing an imperfect bonding parameter on $\p D$. The perfect bonding is characterized by the continuity of the flux and the potential along the interface $\p D$ as given in \eqnref{interface}, while the imperfect bonding is characterized by either discontinuity of the potential or that of the flux along the interface. The former one is referred to as the low conductivity (LC) type, while the latter as the high conductivity (HC) type (see, e.g., \cite{BM}).

The LC type imperfect interface problem is described as follows:
\beq\label{LC}
\begin{cases}
	\nabla\cdot\Gs\nabla u=0 \quad &\mbox{in } D \cup (\Rbb^d \setminus\ol{D}) ,\\
	\Gb(u|_+ - u|_-)= 	\ds \Gs_m \p_\nu u |_+ \quad&\mbox{on } \p D,\\
	\ds \Gs_c \p_\nu u|_- = \Gs_m \p_\nu u|_+ \quad&\mbox{on } \p D,\\
	u(x)-a \cdot x =O(|x|^{-d+1}) \quad&\mbox{as }|x|\rightarrow\infty.
\end{cases}
\eeq
Here, $\Gb$ is the interface parameter of the LC type, which is a non-negative function defined on the interface $\p D$.

\begin{figure}[t!]
\begin{center}
\epsfig{figure=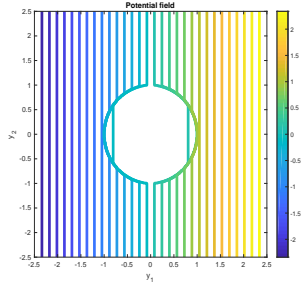, height=4cm, width=5cm}
\end{center}
\caption{Neutral inclusion by imperfect interface \cite{TR}} \label{imperfect}
\end{figure}

It is proved in \cite{TR} that if $D$ is a disk (or a ball) of radius $r$ and
\beq\label{Gbdisk}
\Gb=\frac{1}{r}\frac{\Gs_c\Gs_m}{\Gs_c-\Gs_m},
\eeq
then the solution $u$ to \eqnref{LC} satisfies
\beq\label{neutral2}
u(x)-a \cdot x \equiv 0 \quad\mbox{for all } x \in \Rbb^d \setminus \GO,
\eeq
in other words, $D$ with $\Gb$ is neutral. See Fig. \ref{imperfect}.

It is proved in \cite{KL19}, based on the neutrality criterion obtained in \cite{BM}, that the only neutral inclusions with the imperfect bonding parameters are disks (balls) with constant interface parameters if $\Gs_m$ is isotropic, and ellipses (ellipsoids) if $\Gs_m$ is anisotropic. In the same paper a way to construct an imperfect bonding parameter on the boundary of arbitrary domain has been investigated. For that purpose it is assumed that $D$ is a perfect conductor, meaning that $\Gs_c=\infty$, which is to use the conformal mapping as in \eqnref{f}. Under this assumption, the problem \eqnref{LC} in two dimensions becomes the following one:
\beq\label{weakLC}
\begin{cases}
	\GD u=0 \quad &\mbox{in}~\Rbb^2\setminus\overline{D},\\
	\beta(u-\Gl)= 	\ds \p_\nu u|_{+} \quad&\mbox{on}~\p D,\\
	\vspace{0.3em}
	u(x)-a\cdot x=O(|x|^{-1}) \quad&\mbox{as}~|x|\rightarrow\infty.
\end{cases}
\eeq
The following theorem is obtained.
\begin{thm}[\cite{KL19}]
Let $D$ be a bounded simply connected domain in $\Rbb^2$ with the Lipschitz boundary which admits the conformal mapping $\Phi$ of the form \eqnref{f}. Assume that
\beq\label{bD}
	|b_D|\leq 2-\sqrt{3}.
\eeq
Define $\Gb$ on $\p D$ by
\beq\label{Gb}
\Gb(z) = \left(\frac{1}{1+|b_D|} + \frac{1}{1-|b_D|} -1 + \left(\frac{2}{1+|b_D|} - \frac{2}{1-|b_D|}\right) \cos2\Gt\right) \frac{1}{|\Phi_D'(e^{i\Gt})|}
\eeq
for $z=\Phi_D(e^{i\Gt})$. Then the solution $u$ to the problem \eqnref{weakLC} satisfies $u(x)-a\cdot x=O(|x|^{-2})$ as  $|x|\rightarrow\infty$, namely, $(D,\Gb)$ is weakly neutral to multiple uniform fields.
\end{thm}

\begin{figure}[t!]
\begin{center}
\epsfig{figure=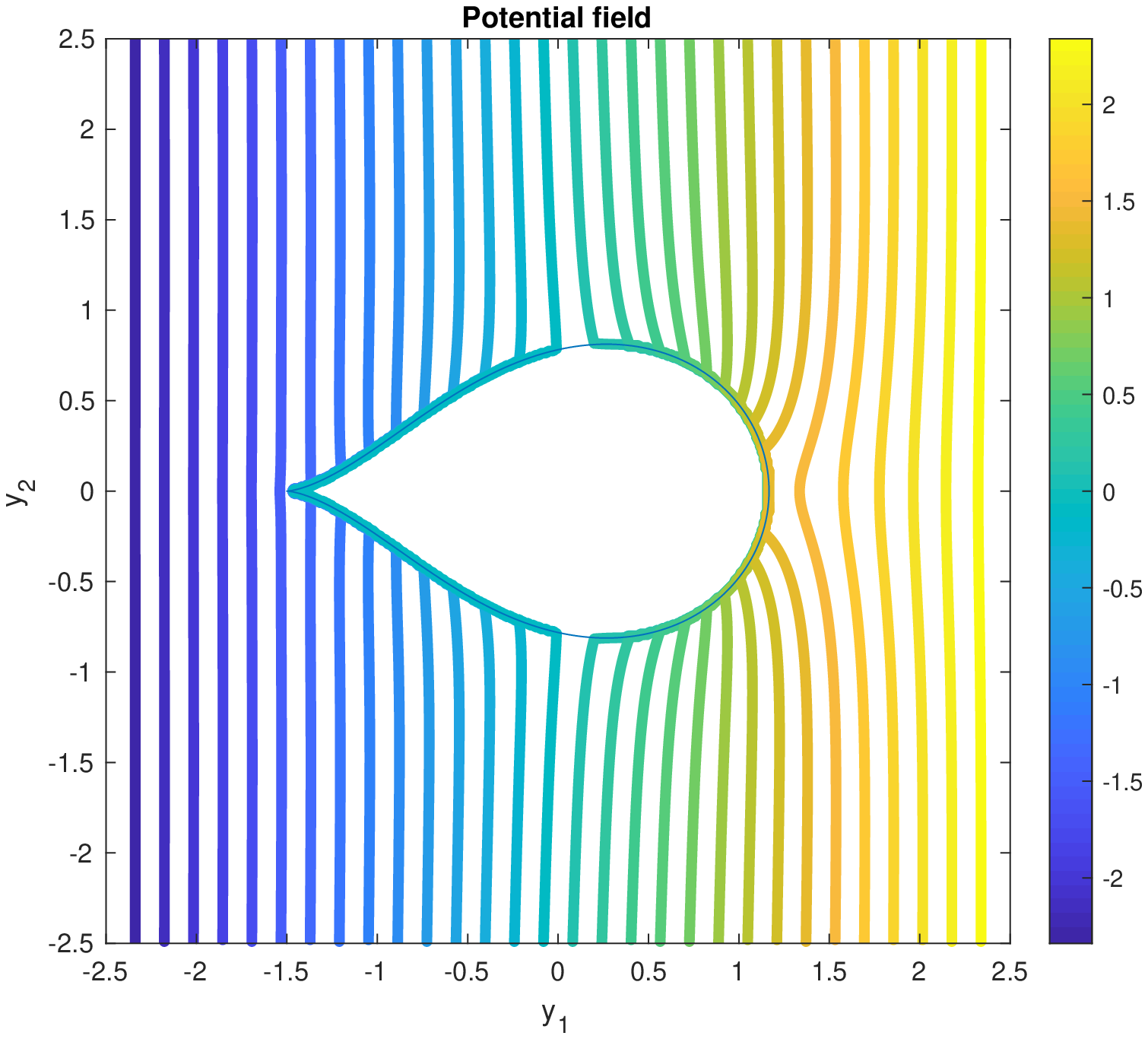, height=3.5cm} \epsfig{figure=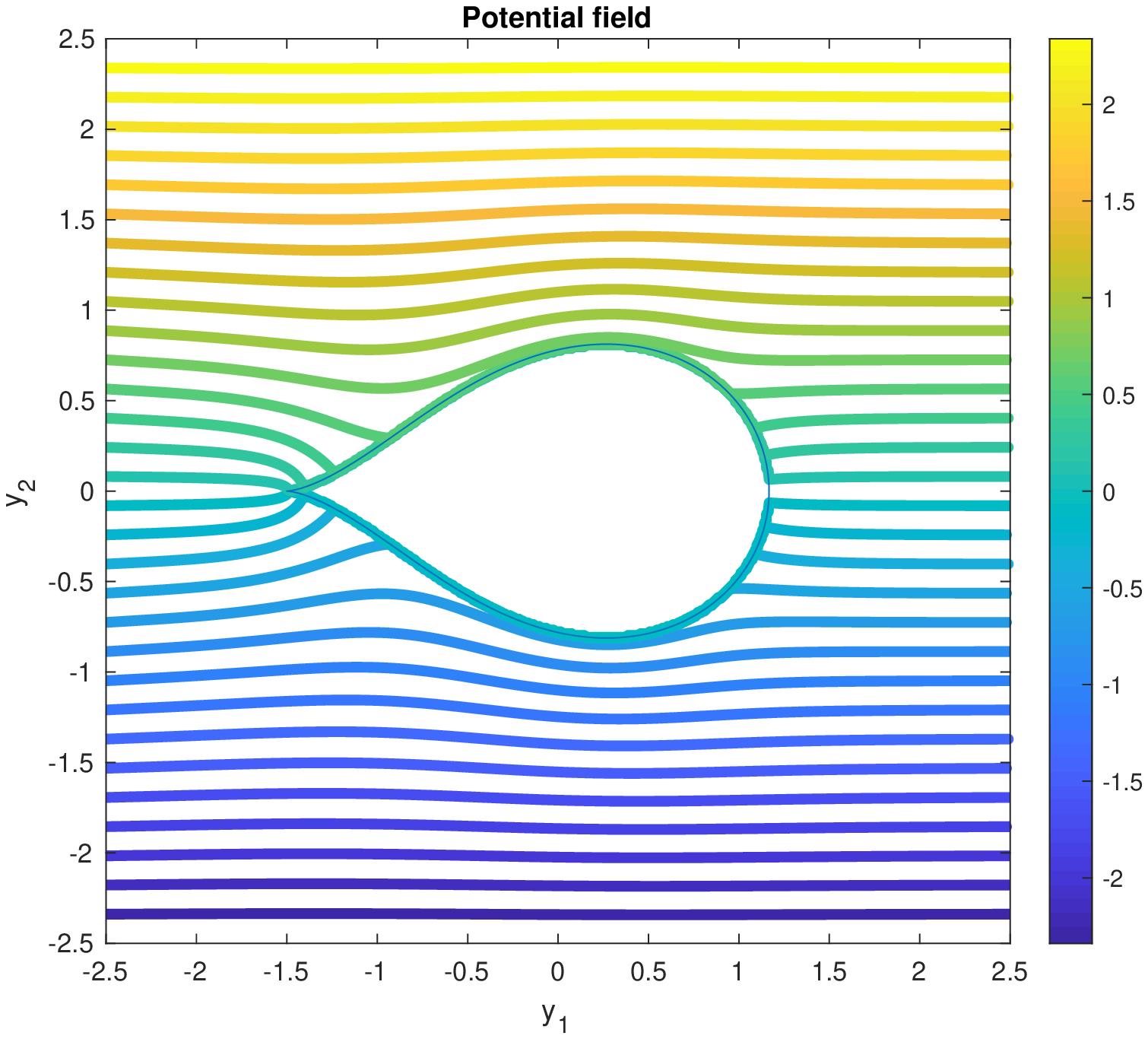, height=3.5cm} \\
\epsfig{figure=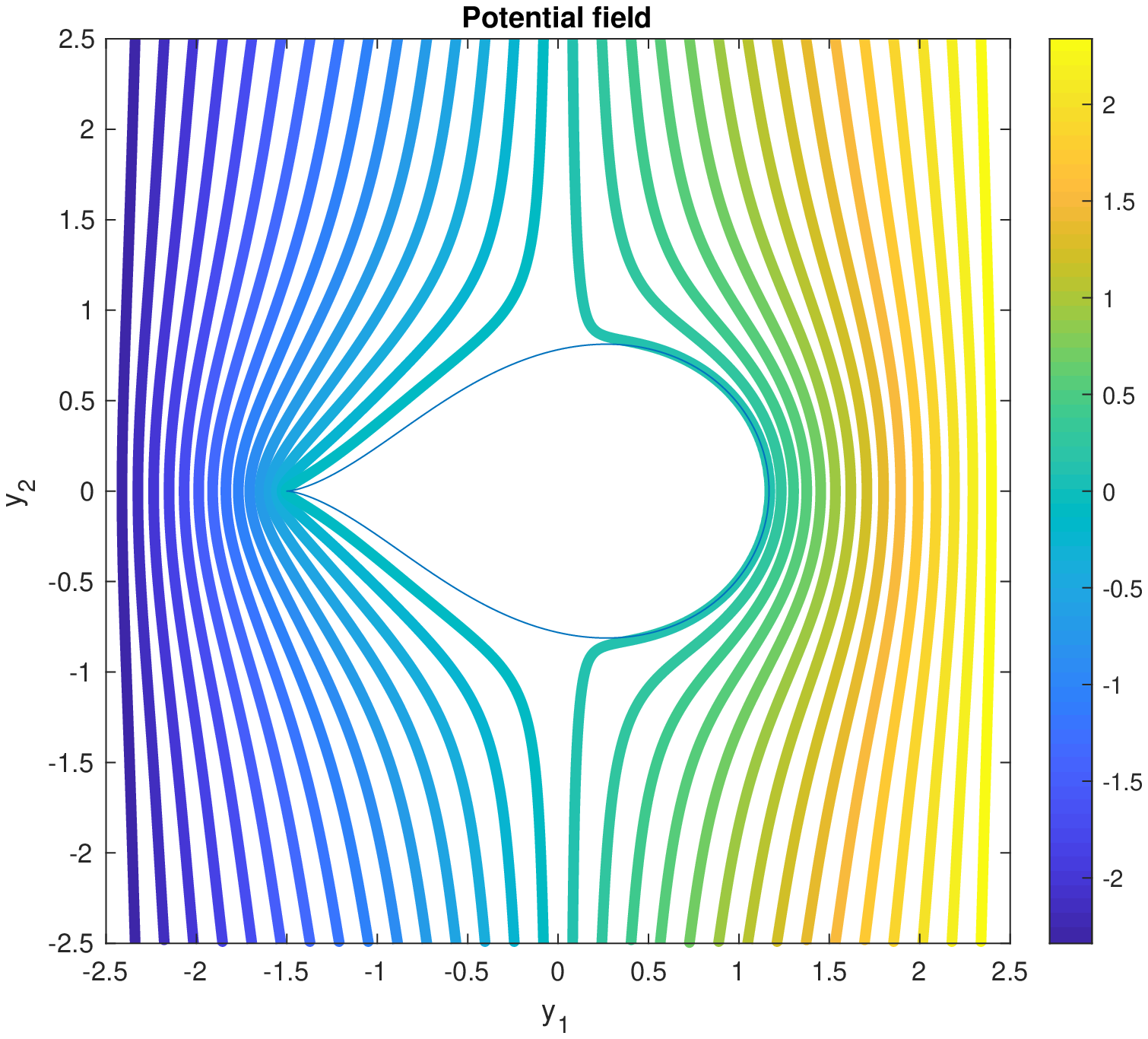, height=3.5cm} \epsfig{figure=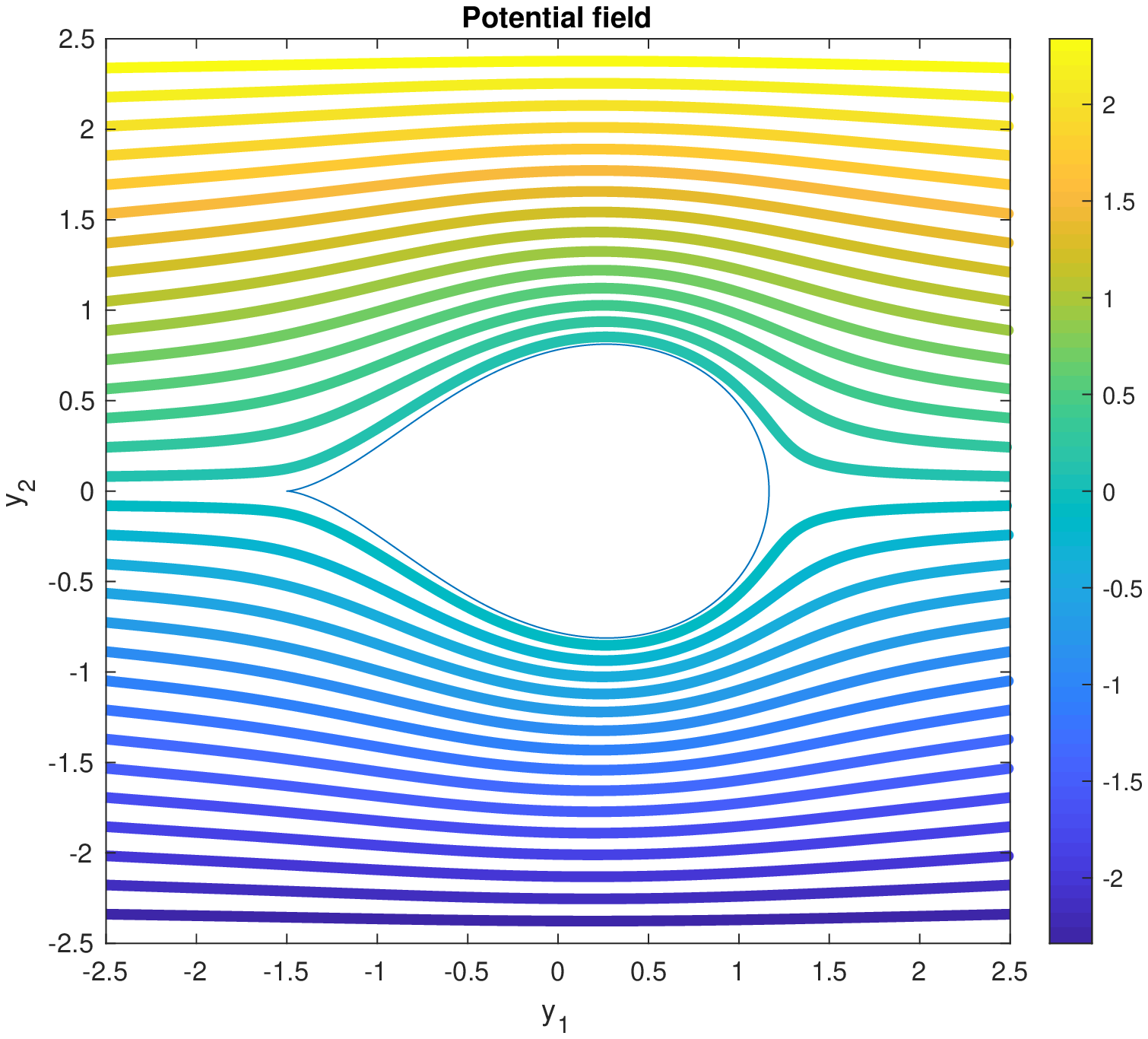, height=3.5cm}
\end{center}
\caption{Upper: the solution with the imperfect bonding parameter, Lower: without it. The solution in the upper one is less perturbed than that in the lower one}\label{imweak}
\end{figure}

It is helpful to mention that the condition \eqnref{bD} is imposed, even though the definition \eqnref{Gb} makes sense without the condition, to guarantee the function $\Gb$ defined by \eqnref{Gb} being positive, and hence uniqueness of the solution to \eqnref{weakLC}. Fig. \ref{imweak} clearly shows that the field with the imperfect bonding parameter is less perturbed than that without it.

\end{document}